\documentclass{article}

\usepackage{amsfonts}
\usepackage{amssymb}
\usepackage{amsthm}
\usepackage{amsmath}
\usepackage{epsfig}
\usepackage{pstricks}
\usepackage{pst-node}
\usepackage{psfrag}
\usepackage{url}
\usepackage{mathdots}
\usepackage{xcolor}
\usepackage{lscape}
\usepackage{rotating}

\newtheorem{theorem}{Theorem}[section]
\newtheorem{corollary}[theorem]{Corollary}
\newtheorem{lemma}[theorem]{Lemma}
\newtheorem{proposition}[theorem]{Proposition}
\newtheorem{properties}[theorem]{Properties}

\newtheorem{definition}[theorem]{Definition}

\def\N{\mathbb{N}}
\def\T{\mathbb{T}}

\def\R{\mathbb{R}}
\def\RR{\mathbb{R}}
\def\real{\mathbb{R}}

\def\Z{\mathbb{Z}}

\def\gauss{\mathcal{G}}
\def\mappair{\mathcal{H}}
\def\conv{\operatorname{conv}}
\def\cone{\operatorname{cone}}

\newcommand{\dist}{\textnormal{dist}}
\newcommand{\st}{\ |\ }
\newcommand{\iso}{\cong}
\newcommand{\wo}{\backslash}
\newcommand{\wt}{\widetilde}

\graphicspath{{Figures/}}

\def \notetoauthors#1{}

\begin{document}

\title{The width of $5$-dimensional prismatoids}

\author{
Benjamin Matschke\thanks{Supported by Deutsche Telekom Stiftung, NSF Grant DMS-0635607, and an EPDI fellowship.
}
\and \setcounter{footnote}{0}
Francisco Santos$^{*}$\thanks{$^{*}$%
Supported in part by the Spanish Ministry of Science through grants  MTM2008-04699-C03-02,  MTM2011-22792 and CSD2006-00032 (i-MATH) and by MICINN-ESF EUROCORES programme EuroGIGA - ComPoSe IP04 - Project EUI-EURC-2011-4306.
}
\and \setcounter{footnote}{0}
Christophe Weibel$^{**}$\thanks{$^{**}$Supported by NSF Grant CCF-1016778.
}}

\date{\today}

\maketitle

\begin{abstract}
Santos' construction of counter-examples to the Hirsch Conjecture (2012) is based on the existence of prismatoids of dimension $d$ of width greater than $d$. Santos, Stephen and Thomas (2012) have shown that this cannot occur in $d\le 4$. Motivated by this we here study the width of $5$-dimensional prismatoids, obtaining the following results:
\begin{itemize}
\item There are $5$-prismatoids of width six with only $25$ vertices, versus the $48$ vertices in Santos' original construction. This leads to non-Hirsch polytopes of dimension $20$, rather than the original dimension $43$.
\item There are $5$-prismatoids with $n$ vertices and width $\Omega(\sqrt{n})$ for arbitrarily large $n$. Hence, the width of $5$-prismatoids is unbounded.
\end{itemize}
\medskip
\noindent\textbf{Mathematics Subject Classification:} 52B05, 52B55, 90C05

\end{abstract}

\section{Introduction}

\paragraph*{The Hirsch Conjecture.}
In 1957, Warren Hirsch conjectured that the (combinatorial) diameter of any $d$-dimensional convex polytope or polyhedron with $n$ facets is at most $n-d$. 
Here the diameter of a polytope is the diameter of its graph, that is, the maximal number of edges that one needs to pass in order to go from one vertex to another.

Hirsch's motivation was that the simplex method, devised by Dantzig ten years earlier, finds the optimal solution of a linear program by walking along the graph of the feasibility polyhedron. In particular, the maximum diameter among polyhedra with a given number of  facets  is a lower bound for the computational complexity of the simplex algorithm for 
\emph{any} pivoting rule. Put differently, our lack of knowledge of whether polynomial-time pivot rules exist for the simplex method (and whether strongly polynomial-time algorithms exist  for linear programming, which is one of Steve Smale's ``mathematical problems for the 21st century''\cite{Smale:problems}) is related to the fact that we do not know whether a polynomial upper bound exists for diameters of polytopes. 
For this reason the Hirsch Conjecture raised a lot of interest both among discrete geometers and optimization theorists. We refer to the comprehensive surveys \cite{Kim-Santos:survey, Klee-Kleinschmidt:survey} and just mention that the best known upper bounds on the diameter of a $d$-polytope with $n$ facets are $n2^{d-2}/3$ (for $n>d\geq 3$) \cite{Larman:diameter} and $n^{\log_2d+1}$ \cite{Kalai-Kleitman}.
The existence of a polynomial upper bound is usually dubbed the \emph{Polynomial Hirsch Conjecture}, and was recently the subject of Polymath project no.~3 under the coordination of Gil Kalai~\cite{Kalai:polymath3}.

\paragraph{Prismatoids and pairs of geodesic maps.}
An unbounded polyhedron violating the Hirsch Conjecture was found in 1967 by Klee and Walkup~\cite{Klee-Walkup:d-step}. But the bounded Hirsch Conjecture resisted until 2010, when Santos constructed a $43$-dimensional polytope with $86$ facets whose diameter was proved to exceed 43~\cite{Santos:Hirsch-counter}.
His construction is based on the following statement.


\begin{theorem}[Santos~\cite{Santos:Hirsch-counter}]
\label{thm:dstep-prismatoid}
If $Q$ is a prismatoid of dimension $d$ with $n$ vertices and width $l$, then there is another prismatoid $Q'$ of dimension $n-d$, with $2n-2d$ vertices and width at least $l+n-2d$.
%
In particular, if $l>d$ then the dual of $Q'$ violates the $d$-step Conjecture, and also the Hirsch Conjecture.
\end{theorem}

Here a \emph{prismatoid} is a polytope having two parallel facets $Q^+$ and $Q^-$ that contain all vertices. We call $Q^+$ and $Q^-$ the \emph{base facets} of $Q$.
The \emph{width} of a prismatoid is the minimum number of steps needed to go from $Q^+$ to $Q^-$, where a step consists in moving from a facet to an adjacent one, crossing a \emph{ridge}.

The only $2$-dimensional prismatoids are the trapezoids, which clearly have width two. Concerning $3$-dimensional prismatoids, it is very easy to prove that they all have width equal to $2$ or $3$, depending on the existence or not of a facet sharing edges with both bases. With a less trivial proof, it is also true that $4$-prismatoids have width at most four~\cite{Santos-Stephen-Thomas:4prismatoid}. So, prismatoids leading to counter-examples to the Hirsch Conjecture via Theorem~\ref{thm:dstep-prismatoid} need to have dimension five or larger. For this reason it seems specially relevant to study the width of $5$-dimensional prismatoids, as we do in this paper.

Our study is based in the following reduction, also from~\cite{Santos:Hirsch-counter}. Let $Q\subset\real^d$ be a prismatoid with bases $Q^+$ and $Q^-$. Consider the bases as simultaneously embedded into $\real^{d-1}$ via an affine projection $\real^d\to \real^{d-1}$ and let $\gauss^+$ and $\gauss^-$ be the intersections of their normal fans with the unit sphere $S^{d-2}$. $\gauss^+$ and $\gauss^-$ are geodesic cell decompositions of $S^{d-2}$, what we call \emph{geodesic maps}. Let $\mappair$ be their common refinement: cells of $\mappair$ are all the intersections of a cell of $\gauss^+$ and a cell of $\gauss^-$. Then:
\begin{itemize}
\item All facets of $Q$ other than the two bases appear as vertices of~$\mappair $.
\item The facets adjacent to $Q^+$ (respectively to $Q^-$) appear in $\mappair $ as the vertices of $\gauss^+$ (respectively of $\gauss^-$).
\item Adjacent facets of $Q$ appear as vertices connected by an edge of~$\mappair $.
\end{itemize}

As a consequence, we have the following result, in which we call \emph{width} of the pair of geodesic maps $(\gauss^+, \gauss^-)$ the minimum graph distance along the graph of~$\mappair$ from a vertex of $\gauss^+$ to a vertex of~$\gauss^-$:

\begin{lemma}[Santos~\cite{Santos:Hirsch-counter}]
\label{lemma:maps}
The width of a prismatoid $Q\subset\real^d$ equals two plus the width of its corresponding pair of maps $(\gauss^+,\gauss^-)$ in $S^{d-2}$.
\end{lemma}

In particular, the width of a $5$-prismatoid~$Q$ is two plus the width of its corresponding pair of maps $(\gauss^+,\gauss^-)$ in the $3$-sphere. The number of vertices of~$Q$ clearly equals the sum of the numbers of maximal cells (called facets) in~$\gauss^+$ and~$\gauss^-$. 

\paragraph*{Main results.}
Our goal in this paper is to construct $5$-prismatoids with large width and few vertices. What we get is:

\begin{itemize}
\item In Section~\ref{sec:FewVertices} we construct several prismatoids  of width six, the smallest ones having $12+13=25$ vertices (Theorem~\ref{thm:q25}). This is much smaller than the original example from~\cite{Santos:Hirsch-counter}, which had $48$ vertices. Applying Theorem~\ref{thm:dstep-prismatoid} to our prismatoid we now have counter-examples to the Hirsch Conjecture in dimension $20$, rather than the original $43$. Moreover, these new non-Hirsch polytopes are small enough to be computed explicitly (see Section~\ref{sub:Explicit}), while the original one was not:

\begin{theorem}
\label{thm:small-non-hirsch-intro}
The $20$-prismatoid whose $40$ vertices are the rows of Table~\ref{table:poly20dim21} has width 21. Hence, its polar is a non-Hirsch $20$-dimensional polytope with $40$ facets. This non-Hirsch polytope has 
$36,425$ vertices.
\end{theorem}

\item In Section~\ref{sec:LargeWidth} we construct an infinite family of $5$-prismatoids of width growing arbitrarily. More precisely:
\begin{theorem}[Corollary~\ref{corPrismatoidWithArbitraryLargeWidth}]
\label{thmPrismatoidWithArbitraryLargeWidth-intro}
For every $k$ there is
a $5$-dimensional prismatoid with $6k(6k-1)$ vertices in each base facet and of width $4+k$.
\end{theorem}
\end{itemize}


Via Theorem~\ref{thm:dstep-prismatoid} these two constructions give counter-examples to the Hirsch Conjecture. Apart of the fact that the first construction yields the smallest non-Hirsch polytope known so far, let us analyze ``how non-Hirsch'' our polytopes are, computing their excess. Following~\cite{Santos:Hirsch-counter} we call \emph{(Hirsch) excess} of a $d$-polytope with $n$ facets and diameter $\delta$ the quantity
\[
\frac{\delta}{n-d} -1.
\]
which is positive if and only if the diameter exceeds the Hirsch bound. Excess is a significant parameter since, as shown in~\cite[Section 6]{Santos:Hirsch-counter}, from any non-Hirsch polytope one can obtain infinite families of them with (essentially) the same excess as the original, even in fixed dimension.

The excess of the non-Hirsch polytope produced via Theorem~\ref{thm:dstep-prismatoid} from a $d$-prismatoid $Q$ of width $l$ and $n$ vertices equals
\[
\frac{l-d}{n-d},
\]
so we call that quotient the \emph{(prismatoid) excess} of $Q$. The prismatoid of Section~\ref{sub:25vertices}, hence also the non-Hirsch polytope of Theorem~\ref{thm:small-non-hirsch-intro}, has excess $1/20$. This is the greatest excess of a prismatoid or polytope constructed so far. (The excess of the Klee-Walkup \emph{unbounded} non-Hirsch polyhedron is, however, $1/4$). 

It could seem that the prismatoids of arbitrarily large width from Theorem~\ref{thmPrismatoidWithArbitraryLargeWidth-intro}
should give greater excess. However, this is not the case. They have excess
\[
\frac{k-1}{12k(6k-1)-5}\ ,
\]
which is much smaller than $1/20$ for small values of $k$ and goes to zero for large values.
In this sense, the second construction does not lead to improved counter-examples to the Hirsch Conjecture, but it is interesting theoretically: It shows a crucial difference between $4$-prismatoids, which have width at most $4$, and $5$-prismatoids, which can have width arbitrarily large.

On the side of upper bounds, Larman's general bound for the diameters of polytopes~\cite{Larman:diameter} implies that the excess of prismatoids of dimension $d$ cannot exceed $2^{d-2} /3$, that is,  $8/3$ in dimension five. In the following result we improve this to $1/3$. 

\begin{proposition}
\label{prop:5prismatoids-bound}
No $5$-prismatoid with $n$ vertices has width larger than $n/3 + 3$.
\end{proposition}

\begin{proof}
Let $n^+$ and $n^-$ denote the numbers of vertices in the two bases of a given prismatoid, and assume without loss of generality that $n^+ \geq n^-$. Let $u$ be a vertex of~$\mappair$ that is not in $\gauss^-$ but at distance $1$ from a vertex of~$\gauss^-$. Then $u$ lies in the intersection of an edge of $\gauss^-$ and a closed $2$-cell $F$ of $\gauss^+$. The rest of the proof concentrates in the polygonal subdivision $\mappair_F$ of $F$ induced by~$\mappair$.

The subdivision $\mappair_F$ has at most $n^-$ $2$-cells. If we look at it as a topological subdivision of a $2$-ball with all vertices of degree at least three (which means that some of the original vertices of the polygon $F$ may not be considered vertices in this topological subdivision) Euler's formula says that it has at most $2n^- - 2$ such vertices. That is to say, $\mappair_F$ has at most $2n^- - 3$ vertices excluding the vertices of the polygon $F$ and the vertex $u$. The key point now is that in $\mappair_F$ there are three disjoint paths from $u$ to vertices of F (see the reason below). The shortest of these three paths uses at most $(2n^- - 3)/3 = (2/3)n^- - 1$ intermediate vertices, hence it has length at most $(2/3)n^-$. This shows that there is a path from a vertex of $\gauss^-$ (a neighbor of $u$) to a vertex of $\gauss^+$ (a vertex of $F$) of length at most $(2/3)n^- + 1\le n/3+1$, which corresponds to a width of at most $n/3+3$.

For the 3 disjoint paths, we argue as follows. $\mappair_F$ is an example of a \emph{regular subdivision} (cf., for example,~\cite{triang-book}). That is, we can lift $F$ to a convex surface in~$\real^3$ whose facets project down to the $2$-cells in which $\gauss^-$ divides $F$. Adding a point ``at infinity" in the projection direction, we get a $3$-polytope whose graph is the graph of $\mappair_F$ with an extra vertex $u_\infty$ joined to the vertices of the original polygon $F$. This graph is $3$-connected, and three disjoint paths from $u$ to $u_\infty$ give what we want.
\end{proof}


Put differently, starting with prismatoids of dimension $5$, Theorem~\ref{thm:dstep-prismatoid} cannot produce polytopes violating the Hirsch Conjecture by more than 33\%.
It would be interesting, however, to know whether there are arbitrarily large $5$-prismatoids with excess bounded away from zero. Put differently, whether constructions similar to those of Section~\ref{sec:LargeWidth} can be done but with the number of vertices growing linearly with respect to the width, rather than quadratically. We believe this to be the case.

\section{Pairs of maps with few vertices} \label{sec:FewVertices}

\subsection{Large width via incidence patterns} \label{sub:IncidencePatterns}

Let us say that a pair of geodesic maps $(\gauss^+,\gauss^-)$ in $S^{d-2}$ has \emph{large width} if its width is larger than $d-2$. For the whole of section~\ref{sec:FewVertices} we assume that~$\gauss^+$ and~$\gauss^-$ are \emph{transversal} in the following sense: for each cell $C^+$ of $\gauss^+$ and cell $C^-$ of $\gauss^-$ with non-empty intersection we have that:
\[
\dim(C^+)  + \dim(C^-) = \dim(C^+ \cap C^-) + (d-2).
\]

The following statement gives a sufficient condition for a transversal
pair of maps to have large width.
For all examples that we know of in the $3$-sphere the condition of the proposition is also necessary,
but we can construct an example in the $4$-sphere where it is not necessary.

\begin{proposition}[Santos~\cite{Santos:Hirsch-counter}]
\label{prop:transversal}
Let $(\gauss^+, \gauss^-)$ be a transversal pair of geodesic maps in the $(d-2)$-sphere. If there is a path of length $d-2$ between a vertex $v_1$ of $\gauss^+$ and a vertex $v_2$ of $\gauss^-$, then the facet of $\gauss^+$ containing $v_1$ in its interior has $v_2$ as a vertex and the  facet of $\gauss^+$ containing $v_2$ in its interior has $v_1$ as a vertex.
\end{proposition}


To use this proposition we introduce the vertex-facet incidence pattern of a pair of maps:

\begin{definition}
\label{defi:incidence}
The \emph{incidence pattern} of a pair of maps $(\gauss^+,\gauss^-)$ in $S^{d-2}$ is the bipartite directed graph having a node for each facet of $\gauss^+$ and $\gauss^-$ and having an arrow from a cell $C$ to a cell $D$ if there is a vertex of $C$ in the interior of $D$.
The \emph{reduced incidence pattern} is the subgraph induced by facets of one map that contain some vertex of the other map.
\end{definition}

The graph is bipartite since all arrows go from a facet of $\gauss^+$ to one of $\gauss^-$ or vice-versa. With this notion, Proposition~\ref{prop:transversal} becomes:

\begin{proposition}
\label{prop:incidence}
Let $(\gauss^+, \gauss^-)$ be a transversal pair of geodesic maps in the $(d-2)$-sphere. If there is no directed cycle of length two in its reduced incidence pattern then the pair has large width.
\end{proposition}

\begin{proof}
For the total (i.e., not reduced) incidence pattern this statement is an exact translation of Proposition~\ref{prop:transversal}.
The reason why we can use the \emph{reduced} incidence pattern in the statement is that
facets of $\gauss^+$ or $\gauss^-$ that do not contain vertices of the other map correspond to sources in the incidence pattern, which cannot participate in directed cycles.
\end{proof}

The original pair of geodesic maps  $(\gauss^+,\gauss^-)$ of large width constructed in~\cite{Santos:Hirsch-counter} has the following reduced incidence pattern:

\begin{itemize}
\item $\gauss^+$ has four facets $A_1$, $A_2$, $C_{12}$ and $C_{34}$ containing all the vertices of $\gauss^-$.
\item $\gauss^-$ has four facets $B_1$, $B_2$, $D_1$ and $D_2$ containing all the vertices of $\gauss^+$.
\item The reduced incidence pattern induced has the sixteen arrows of the form $A_i\to B_j$, $B_i\to C_j$, $C_i\to D_j$,  and $D_i\to A_j$, with $i,j\in\{1,2\}$.
\end{itemize}

This is exactly the pattern that we will reproduce in all the examples in this section. The reason for using this pattern is two-fold. On the one hand, its large symmetry allows us to use symmetry in the construction as well. On the other hand, the pattern is minimal in the following sense. Observe that every vertex in a reduced incidence pattern has out-degree at least two.

\begin{proposition}
\label{prop:minimal}
Let $G$ be a directed bipartite graph with no cycle of length two and with out-degree at least two in every vertex. Then $G$ has at least eight vertices and if it has eight vertices then it has four on each side.
\end{proposition}

\begin{proof}
We denote by $X$ and $Y$ the two parts (subsets of vertices) of $G$.
First note that $|X|\ge 3$ and $|Y|\ge 3$. Indeed, if we have an arrow from an $x\in X$ to a certain $y\in Y$ then the (at least two) arrows coming back from $y$ must go to two vertices of $X$ different from $x$. So, in order to prove that $|X \cup Y| \ge 8$ we only need to consider the case where one of the parts, say $X$, has only three vertices. If this happens, then each vertex of $Y$ has in-degree at most equal to one, since they have out-degree at least two. Since there are at least six arrows from $X$ to $Y$, we have $|Y|\ge 6$ and $|X \cup Y| \ge 9$.
\end{proof}

\begin{figure}
\centerline{
\input{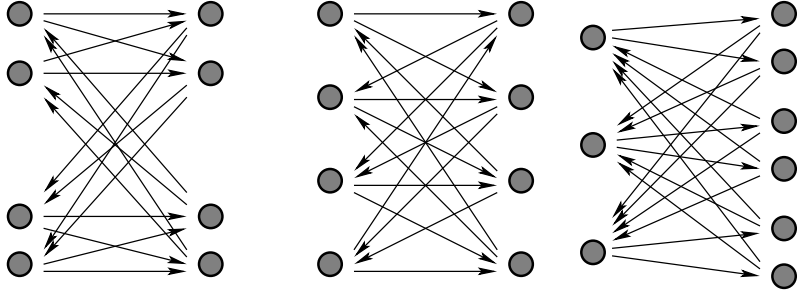}
}
\caption{The two reduced incidence patterns without $2$-cycles and with $4+4$ facets, and the one with $3+6$ facets}
\label{fig:patterns}
\end{figure}

The proof shows that in fact the reduced incidence pattern with $4+4$ facets that we describe above is almost unique. If a directed bipartite graph with eight vertices has no $2$-cycle and minimum out-degree at least two, then its average in-degree is also two, which implies every vertex has in-degree and out-degree \emph{exactly} two. That is, the digraph is obtained by giving directions to the edges of $K_{4,4}$ in such a way that every vertex has in-degree and out-degree equal to two. Equivalently, by decomposing the edge set of $K_{4,4}$ into two $2$-regular subgraphs. There are only two ways of doing that: either each subgraph is two disjoint cycles of length four, or each subgraph is a cycle of length eight (see Figure~\ref{fig:patterns}). The former corresponds to the reduced incidence pattern we described.

Similar arguments show that there is a unique reduced incidence pattern that is possible with three facets in one map and six in the other. It is obtained orienting the edges of the complete $K_{3,6}$ as follows. From each vertex on the small part we have two arrows to vertices in the big part, with no vertex in the big part getting two arrows.

\subsection{Constructions based on the $4$-cube} \label{sub:4cube}

\subsubsection*{A construction with $40$ vertices}

To guarantee that our geodesic maps $\gauss^+$ and $\gauss^-$ in $S^3$ are polytopal we actually construct them as the \emph{central fans} of certain $4$-polytopes $P^+$ and $P^-$. Let:
\begin{eqnarray}
P^+ &:=&\conv\left\{(\pm a, \pm b, \pm c, \pm d), (\pm e, \pm f, \pm g, \pm h) \right \}, \label{eq:P^+}
\\
P^- &:=&\conv\left\{(\pm d, \pm c, \pm a, \pm b), (\pm h, \pm g, \pm e, \pm f) \right \},  \label{eq:P^-}
\end{eqnarray}
$P^+$ and $P^-$ are congruent  via the map $(x_1,x_2,x_3,x_4)$ $\mapsto $   $(x_4,x_3,x_1,x_2)$, and each of them is the common convex hull of two $4$-cubes. We pose the following restrictions on the parameters $a,b,c,d,e,f,g,h\in(0,\infty)$:
\[
a> e, \qquad b> f, \qquad c> g, \qquad d< h.
\]
These restrictions are enough to determine the combinatorial type of~$P^+$ and~$P^-$. Indeed, $P^+$ has the following $20$ facets:
\begin{itemize}
\item Six of the facets of the $4$-cube $\conv\left\{(\pm a, \pm b, \pm c, \pm d)\right\}$, those  given by the inequalities
\begin{equation}
\pm x_1 \le a,\qquad \pm x_2 \le b,\qquad \pm x_3 \le c.
\label{eq:F1-F2-F3}
\end{equation}
We denote them $F_1$, $F_{\overline 1}$, $F_2$, $F_{\overline 2}$, $F_3$, and $F_{\overline 3}$. The index denotes which coordinate remains constant in the facet, while the bar or the absence of it indicates whether the coordinate is negative or positive.

\item Two of the facets of the $4$-cube $\conv\left\{(\pm e, \pm f, \pm g, \pm h)\right\}$, those  given by the inequalities
\begin{equation}
\pm x_4 \le h.
\label{eq:F4}
\end{equation}
We denote them $F_4$, $F_{\overline 4}$.
\item Six times two facets connecting the first six to the last two, given by the inequalities
\begin{eqnarray}
\pm (h-d)  x_1 \pm (a-e) x_4 &\le& ah -ed,  \label{eq:F14}
\\
\pm (h-d)  x_2 \pm (b-f) x_4 &\le& bh -fd,  \label{eq:F24}
\\
\pm (h-d)  x_3 \pm (c-g) x_4 &\le& ch -gd.  \label{eq:F34}
\end{eqnarray}
Following similar conventions, we denote these twelve facets as:
\[
F_{14}, F_{\overline 1 4}, F_{1 \overline 4}, F_{\overline 1 \overline 4}, \qquad
F_{24}, F_{\overline 2 4}, F_{2 \overline 4}, F_{\overline 2 \overline 4}, \qquad
F_{34}, F_{\overline 3 4}, F_{3 \overline 4}, F_{\overline 3 \overline 4}.
\]
For example, $F_{1\overline 4}=\conv\{(a, \pm b, \pm c, -d),(e, \pm f, \pm g, -h)\}$.

\end{itemize}
By symmetry, we get that $P^-$ has also 20 facets, that we denote with similar conventions:
\begin{eqnarray*}
&F'_{1}, F'_{\overline 1}, F'_{2}, F'_{\overline  2}, F'_{3}, F'_{\overline 3}, F'_{4}, F'_{\overline  4},&\\
F'_{12}, F'_{\overline 1 2}, F'_{1 \overline 2}, F'_{\overline 1 \overline 2},& \quad
F'_{13}, F'_{\overline 1 3}, F'_{1 \overline 3}, F'_{\overline 1 \overline 3}, \quad&
F'_{14}, F'_{\overline 1 4}, F'_{1 \overline 4}, F'_{\overline 1 \overline 4}.
\end{eqnarray*}

It is interesting to observe that all the facets are combinatorial $3$-cubes (that is, $P^+$ and $P^-$ are \emph{cubical polytopes}). This property was already present in the construction in~\cite{Santos:Hirsch-counter}, which was designed with similar ideas.%
\footnote{The construction of~\cite{Santos:Hirsch-counter} is essentially the same one as here except there $a> e$, $b< f$, $c> g$, $d< h$.}
In what follows we consider the central projections of the face lattices of $P^+$ and $P^-$ to the unit sphere $S^3$. That is, for each facet $F$ of either $P^+$ or $P^-$ we consider the cone $\cone(F):=\{\lambda p : p\in F, \lambda \in [0,\infty)\}$ and its intersection to the unit sphere. This gives us two maps $\gauss^+$ and $\gauss^-$ in $S^3$,  the central fans of $P^+$ and $P^-$.

We are interested in the reduced incidence pattern of this pair of maps. That is, which facets of one map contain which vertices of the other or equivalently, which cones $\cone(F)$ of facets of $P^+$ contain which vertices of $P^-$, and vice-versa.

\begin{theorem}
\label{thm:40vertices}
Assume that $a,b,c,d,e,f,g,h\in(0,\infty)$ satisfy the following inequalities:
\begin{eqnarray}
a> e, && b> f, \qquad c> g, \qquad d< h \label{eq:40vertices}
\\ \frac{e}{c} &>& \max\left\{ \frac{h}{a}, \frac{g}{b}, \frac{f}{d} \right \} \label{eq:F'1-F3}
\\ \frac{b}{c} &>& \max\left\{ \frac{a}{d}, \frac{c}{a}, \frac{d}{b} \right \} \label{eq:F3-F'2}
\\ \frac{b}{h} &>& \max\left\{ \frac{d}{e}, \frac{c}{f}, \frac{a}{g} \right \} \label{eq:F'2-F4}
\\ \frac{e}{h} &>& \max\left\{ \frac{f}{g}, \frac{g}{e}, \frac{h}{f} \right \} \label{eq:F4-F'1}
\end{eqnarray}
Then, the reduced incidence pattern of the pair $(\gauss^+, \gauss^-)$ is:
\bigskip

\centerline{
\begin{picture}(0,0)%
\includegraphics{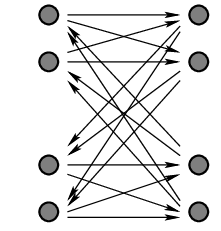}%
\end{picture}%
\setlength{\unitlength}{3947sp}%
\begingroup\makeatletter\ifx\SetFigFont\undefined%
\gdef\SetFigFont#1#2#3#4#5{%
  \reset@font\fontsize{#1}{#2pt}%
  \fontfamily{#3}\fontseries{#4}\fontshape{#5}%
  \selectfont}%
\fi\endgroup%
\begin{picture}(1755,1836)(1561,-3289)
\put(3301,-2836){\makebox(0,0)[lb]{\smash{{\SetFigFont{12}{14.4}{\rmdefault}{\mddefault}{\updefault}{\color[rgb]{0,0,0}$F_4$}%
}}}}
\put(3301,-1636){\makebox(0,0)[lb]{\smash{{\SetFigFont{12}{14.4}{\rmdefault}{\mddefault}{\updefault}{\color[rgb]{0,0,0}$F_3$}%
}}}}
\put(3301,-3211){\makebox(0,0)[lb]{\smash{{\SetFigFont{12}{14.4}{\rmdefault}{\mddefault}{\updefault}{\color[rgb]{0,0,0}$F_{\overline 4}$}%
}}}}
\put(1576,-3211){\makebox(0,0)[lb]{\smash{{\SetFigFont{12}{14.4}{\rmdefault}{\mddefault}{\updefault}{\color[rgb]{0,0,0}$F'_{\overline 2}$}%
}}}}
\put(1576,-2011){\makebox(0,0)[lb]{\smash{{\SetFigFont{12}{14.4}{\rmdefault}{\mddefault}{\updefault}{\color[rgb]{0,0,0}$F'_{\overline 1}$}%
}}}}
\put(3301,-2011){\makebox(0,0)[lb]{\smash{{\SetFigFont{12}{14.4}{\rmdefault}{\mddefault}{\updefault}{\color[rgb]{0,0,0}$F_{\overline 3}$}%
}}}}
\put(1576,-2836){\makebox(0,0)[lb]{\smash{{\SetFigFont{12}{14.4}{\rmdefault}{\mddefault}{\updefault}{\color[rgb]{0,0,0}$F'_2$}%
}}}}
\put(1576,-1636){\makebox(0,0)[lb]{\smash{{\SetFigFont{12}{14.4}{\rmdefault}{\mddefault}{\updefault}{\color[rgb]{0,0,0}$F'_1$}%
}}}}
\end{picture}%

}
\end{theorem}

\begin{proof}
The inequalities~(\ref{eq:40vertices}) guarantee that the face lattices of $\gauss^+$ and $\gauss^-$ are as described above. For the incidence pattern, let us  first consider the cones over the facets $F_3$ and $F_{\overline 3}$. From the description of $P^+$ we have that
\[
\cone(F_3) :=\cone\{(\pm a, \pm b, c, \pm d)\}=\left\{ \frac{x_3}{c} \ge \max \left\{\frac{\pm x_1}{a}, \frac{\pm x_2}{b}, \frac{\pm x_4}{d} \right\}\right\}.
\]
Similarly,
\[
\cone(F_{\overline 3}) :=\cone\{(\pm a, \pm b, -c, \pm d)\}=\left\{ -\frac{x_3}{c} \ge \max \left\{\frac{\pm x_1}{a}, \frac{\pm x_2}{b}, \frac{\pm x_4}{d} \right\}\right\}.
\]
Then, inequalities~(\ref{eq:F'1-F3}) guarantee that the vertices of $F'_1:=\conv\{(h,\pm g,\pm e,\pm f)\}$ and $F'_{\overline 1}:=\conv\{(-h,\pm g,\pm e,\pm f)\}$ are all in the cones $\cone(F_3)$ or $\cone(F_{\overline 3})$.

With the same arguments:
Inequalities~(\ref{eq:F3-F'2}) guarantee that the vertices of $F_3:=\conv\{(\pm a,\pm b, c, \pm d)\}$ and $F_{\overline 3}:=\conv\{(\pm a,\pm b,-c,\pm d)\}$ are all in the cones
\[
\cone(F'_2) :=\cone\{(\pm d, c, \pm a, \pm b)\}=\left\{ \frac{x_2}{c} \ge \max \left\{\frac{\pm x_1}{d}, \frac{\pm x_3}{a}, \frac{\pm x_4}{b} \right\}\right\}
\]
or
\[
\cone(F'_{\overline 2}) :=\cone\{(\pm d,- c, \pm a, \pm b)\}=\left\{ -\frac{x_2}{c} \ge \max \left\{\frac{\pm x_1}{d}, \frac{\pm x_3}{a}, \frac{\pm x_4}{b} \right\}\right\};
\]
Inequalities~(\ref{eq:F'2-F4}) guarantee that the vertices of $F'_2:=\conv\{(\pm d, c, \pm a, b)\}$ and $F'_{\overline 2}:=\conv\{(\pm d, -c,\pm a,\pm b)\}$ are all in the cones
\[
\cone(F_4) :=\cone\{(\pm e, \pm f, \pm g,  h)\}=\left\{ \frac{x_4}{h} \ge \max \left\{\frac{\pm x_1}{e}, \frac{\pm x_2}{f}, \frac{\pm x_3}{g} \right\}\right\}
\]
or
\[
\cone(F_{\overline 4}) :=\cone\{(\pm e,\pm f, \pm g, -h)\}=\left\{ -\frac{x_4}{h} \ge \max \left\{\frac{\pm x_1}{e}, \frac{\pm x_2}{f}, \frac{\pm x_3}{g} \right\}\right\};
\]
And
inequalities~(\ref{eq:F4-F'1}) guarantee that the vertices of $F_4:=\conv\{(\pm e,\pm f, \pm g, h)\}$ and $F_{\overline 4}:=\conv\{(\pm e,\pm f,\pm g, -h)\}$ are all in the cones
\[
\cone(F'_1) :=\cone\{(h, \pm g, \pm e, \pm f)\}=\left\{ \frac{x_1}{h} \ge \max \left\{\frac{\pm x_2}{g}, \frac{\pm x_3}{e}, \frac{\pm x_4}{f} \right\}\right\}
\]
or
\[
\cone(F'_{\overline 1}) :=\cone\{(- h,\pm g, \pm e, \pm f)\}=\left\{ -\frac{x_1}{h} \ge \max \left\{\frac{\pm x_2}{g}, \frac{\pm x_3}{e}, \frac{\pm x_4}{f} \right\}\right\}.
\]
\end{proof}

The following are values of the eight parameters for which all the inequalities are satisfied:
\[
a=6, \quad b=10, \quad c=3, \quad d=2, \quad
e=5, \quad f=3, \quad g=2, \quad h=3.
\]
Hence:

\begin{corollary}
\label{coro:40vertices}
Let
$P^+=\conv\{(\pm6,\pm10,\pm3,\pm2), (\pm5,\pm3,\pm2,\pm3)\}$ and let
$P^-=\conv\{(\pm3,\pm2,\pm5,\pm3), (\pm2,\pm3,\pm6,\pm10)\}$. Let $Q^+$ and $Q^-$ be the polars of $P^+$ and $P^-$ and let
\[
Q_{40}:=\conv (\,Q^+\times \{1\},\ Q^-\times\{-1\}\,).
\]
Then, $Q_{40}$ is a $5$-prismatoid with $40$ vertices and of width (at least) six.
\end{corollary}

The reason why we need the polars of $P^+$ and $P^-$ is that the central fan of a polytope equals the normal fan of its polar.
It can be checked computationally that the width of $Q_{40}$ is \emph{exactly} six.

\subsubsection*{From $40$ to $28$ vertices} \label{sub:28vertices}

In order to reduce the number of vertices of our prismatoid, that is, the number of facets of $P^+$ and $P^-$, we play with the parameters $a,b,c,d,e,f,g,h$. Observe that if we weaken the strict inequalities $a>e$, $b>f$, $c>g$ and $d<h$ to be non-strict inequalities then the combinatorics of $P^+$ and $P^-$ changes only in the direction of merging some of the facets, hence reducing the number of them. In particular, if we let
\[
a= e, \qquad b> f, \qquad c> g, \qquad d< h
\]
the (only) changes that we get are:
\begin{itemize}
\item The three facets $F_1$, $F_{14}$ and $F_{1\overline4}$ of $P^+$ become a single facet still corresponding to the inequality $x_1\le a$
(equations~(\ref{eq:F14}) are now redundant)
and with vertex set $\{(a,\pm b, \pm c, \pm d)$, $(a, \pm f, \pm g, \pm h)\}$. We still denote this facet $F_1$.
\item Similarly, the three facets $F_{\overline 1}$, $F_{\overline 14}$ and $F_{\overline 1\overline4}$ of $P^+$ become a single facet $F_{\overline 1}$.
\item The same changes occur in $P^-$: the three facets $F'_3$, $F'_{13}$ and $F_{\overline13}$ merge into a facet that we still denote
$F'_3$ and the three facets $F'_{\overline 3}$, $F'_{1 \overline 3}$ and $F_{\overline 1 \overline 3}$ merge into a facet that we still denote
$F'_{\overline 3}$.
\end{itemize}

Hence, we now have polytopes $P^+$ and $P^-$ with $16$ facets each. Although we do not need this property, observe that $P^+$ is now a prism with bases $F_1$ and $F_{\overline 1}$, and  $P^-$ is a prism with bases $F_3$ and $F_{\overline 3}$.

 The good news is that (the descriptions of) the cones over the facets do not change at all, except that the cones of merged facets merge. But none of the facets that are merged were involved in the reduced incidence pattern that we proved in Theorem~\ref{thm:40vertices} so we automatically get:

\begin{theorem}
\label{thm:32vertices}
Assume that $a,b,c,d,e,f,g,h$ satisfy all the inequalities of Theorem~\ref{thm:40vertices} except that $a>e$ is changed to $a=e$. Then the reduced incidence pattern is the same as in Theorem~\ref{thm:40vertices}.
\qed
\end{theorem}

\begin{corollary}
\label{coro:32vertices}
Let
$P^+=\conv\{(\pm5,\pm8,\pm3,\pm2), (\pm5,\pm3,\pm2,\pm3)\}$ and let
$P^-=\conv\{(\pm3,\pm2,\pm5,\pm3), (\pm2,\pm3,\pm5,\pm8)\}$. Let $Q^+$ and $Q^-$ be the polars of $P^+$ and $P^-$ and let
\[
Q_{32}:=\conv (\,Q^+\times \{1\}, Q^-\times\{-1\}\,).
\]
Then, $Q_{32}$ is a $5$-prismatoid with $32$ vertices and of width at least six.
\end{corollary}

\begin{proof}
Check that the values
\[
a=5, \quad b=8, \quad c=3, \quad d=2, \quad
e=5, \quad f=3, \quad g=2, \quad h=3
\]
satisfy all the inequalities of Theorem~\ref{thm:40vertices} except that $a=e$ instead of $a>e$.
\end{proof}

Before going further, let us compute the explicit description of $Q_{32}$. The facet descriptions of $P^+$ and $P^-$ are inequalities~(\ref{eq:F1-F2-F3}) to~(\ref{eq:F34}), except that~(\ref{eq:F14}) are redundant when $a=e$. That is:
\[
P^+ =
\left\{
\begin{array}{rcl}
\pm x_1 &\le& 5 \\
\pm x_2 &\le& 8 \\
\pm x_3 &\le& 3 \\
\pm x_4 &\le& 3 \\
\pm   x_2 \pm 5 x_4 &\le& 18 \\
\pm   x_3 \pm  x_4 &\le& 5 \\
\end{array}
\right\}
=
\left\{
\begin{array}{rcl}
\pm 72 x_1 &\le& 360 \\
\pm 45 x_2 &\le& 360 \\
\pm 120 x_3 &\le& 360 \\
\pm 120 x_4 &\le& 360 \\
\pm  20 x_2 \pm 100 x_4 &\le& 360 \\
\pm   72 x_3 \pm  72 x_4 &\le& 360 \\
\end{array}
\right\}.
\]
Similarly,
\[
P^- =
\left\{
\begin{array}{rcl}
\pm x_3 &\le& 5 \\
\pm x_4 &\le& 8 \\
\pm x_2 &\le& 3 \\
\pm x_1 &\le& 3 \\
\pm   x_4 \pm 5 x_1 &\le& 18 \\
\pm   x_2 \pm  x_1 &\le& 5 \\
\end{array}
\right\}
=
\left\{
\begin{array}{rcl}
\pm 72 x_3 &\le& 360 \\
\pm 45 x_4 &\le& 360 \\
\pm 120 x_2 &\le& 360 \\
\pm 120 x_1 &\le& 360 \\
\pm  20 x_4 \pm 100 x_1 &\le& 360 \\
\pm   72 x_2 \pm  72 x_1 &\le& 360 \\
\end{array}
\right\}.
\]
We have normalized the right hand sides so that the polar polytopes~$Q^+$ and~$Q^-$ are (modulo a global scaling factor) the convex hulls of the left-hand side coefficient vectors. That is:
{\footnotesize
\[
Q_{32}:= \operatorname{conv}
\left\{
\bordermatrix{
&x_1&x_2&x_3&x_4&x_5\cr
&  \pm72&   0&   0&   0&   1 \cr
&  0&   \pm45&   0&   0&   1 \cr
&   0&   0&  \pm120&   0&   1 \cr
&   0&   0&   0&  \pm120&   1 \cr
&  0&  \pm 20&  0&   \pm 100&   1 \cr
&   0&   0&  \pm72& \pm 72&   1 \cr
}
,\
\bordermatrix{
&x_1&x_2&x_3&x_4&x_5\cr
&   0&   0&  \pm72&   0&   -1 \cr
&   0&   0&  0&   \pm45&   -1 \cr
&   0&   \pm120&   0&   0&   -1 \cr
&   \pm120&   0&   0&   0&   -1 \cr
&  \pm 100&  0&   0&   \pm 20 &   -1 \cr
&\pm72& \pm72&   0&   0&     -1 \cr
}
\right\}.
\]
}

To go down to $28$ vertices we notice that the polytopes $P^+$ and $P^-$ still have some superfluous facets. In fact, another way of interpreting equations~(\ref{eq:F3-F'2}) and~(\ref{eq:F'2-F4}) of Theorem~\ref{thm:40vertices} is that they show, respectively:
\begin{eqnarray*}
\cone(F_2)\subseteq \cone(F'_2),&&
\cone(F_{\overline 2})\subseteq \cone(F'_{\overline 2}),\\
\cone(F'_4)\subseteq \cone(F_4),&&
\cone(F'_{\overline 4})\subseteq \cone(F_{\overline 4}).
\end{eqnarray*}
This implies that if we remove the inequalities corresponding to the four facets
$F_2$, $F_{\overline 2}$, $F'_4$ and $F'_{\overline 4}$ from the facet definitions of $P^+$ and $P^-$,
the reduced incidence pattern does not change. That is, we now consider
\[
P^+ =
\left\{
\begin{array}{rcl}
\pm x_1 &\le& 5 \\
\pm x_3 &\le& 3 \\
\pm x_4 &\le& 3 \\
\pm   x_2 \pm 5 x_4 &\le& 18 \\
\pm   x_3 \pm  x_4 &\le& 5 \\
\end{array}
\right\}
=
\left\{
\begin{array}{rcl}
\pm 18 x_1 &\le& 90 \\
\pm 30 x_3 &\le& 90 \\
\pm 30 x_4 &\le& 90 \\
\pm  5 x_2 \pm 25 x_4 &\le& 90 \\
\pm   18 x_3 \pm  18 x_4 &\le& 90 \\
\end{array}
\right\}.
\]
and
\[
P^- =
\left\{
\begin{array}{rcl}
\pm x_3 &\le& 5 \\
\pm x_2 &\le& 3 \\
\pm x_1 &\le& 3 \\
\pm   x_4 \pm 5 x_1 &\le& 18 \\
\pm   x_2 \pm  x_1 &\le& 5 \\
\end{array}
\right\}
=
\left\{
\begin{array}{rcl}
\pm 18 x_3 &\le& 90 \\
\pm 30 x_2 &\le& 90 \\
\pm 30 x_1 &\le& 90 \\
\pm  5 x_4 \pm 25 x_1 &\le& 90 \\
\pm   18 x_2 \pm  18 x_1 &\le& 90 \\
\end{array}
\right\}.
\]

$P^+$ and $P^-$ are still prisms over certain $3$-polytopes with $12$ facets.
Their vertex descriptions are:
\begin{eqnarray*}
P^+=\conv\{(\pm5,\pm8,\pm3,\pm2), (\pm5,\pm3,\pm2,\pm3), (\pm5,\pm18,\pm3,0)\},\\
P^-=\conv\{(\pm3,\pm2,\pm5,\pm3), (\pm2,\pm3,\pm5,\pm8), (0,\pm3,\pm5,\pm18)\}.
\end{eqnarray*}

Going to the polars we get:

\begin{corollary}
\label{coro:28vertices}
The following $5$-prismatoid, with $28$ vertices, has width (at least) six:
{\footnotesize
\[
Q_{28}:= \operatorname{conv}
\left\{
\bordermatrix{
&x_1&x_2&x_3&x_4&x_5\cr
&  \pm18&   0&   0&   0&   1 \cr
&   0&   0&  \pm30&   0&   1 \cr
&   0&   0&   0&  \pm30&   1 \cr
&  0&  \pm 5&  0&   \pm 25&   1 \cr
&   0&   0&  \pm18& \pm 18&   1 \cr
}
\quad
\bordermatrix{
&x_1&x_2&x_3&x_4&x_5\cr
&   0&   0&  \pm18&   0&   -1 \cr
&   0&   \pm30&   0&   0&   -1 \cr
&   \pm30&   0&   0&   0&   -1 \cr
&  \pm 25&  0&   0&   \pm 5 &   -1 \cr
&\pm18& \pm18&   0&   0&     -1 \cr
}
\right\}.
\]
}
\end{corollary}

\subsection{Further reduction of the number of vertices} \label{sub:25vertices}

To reduce further the number of vertices of the prismatoid, we are
forced to abandon symmetry. We reduce the number of facets of the
polytopes $P^+$ and $P^-$ by merging some of their facets, while
maintaining the same reduced incidence pattern. In order to keep
control of the reduced incidence pattern, we identify some desirable
properties.

For any polytope $Q$ in $\real^4$, we denote as $\pi_{12}(Q)$ its
orthogonal projection on the plane $P_{12}:=\{x_3=x_4=0\}$, and respectively
denote as $\pi_{34}(Q)$ its orthogonal projection on the plane
$P_{34}:=\{x_1=x_2=0\}$. The projection of $Q$ on a plane is a
polygon whose vertices are projections of some of the vertices
of~$Q$.  We define $C_{12}:= S^1\times
\{(0,0)\}= S^3\cap P_{12}$ and $C_{34}:=\{(0,0)\}\times S^1= S^3\cap
P_{34}$. We also denote as $\gauss(Q)$ the geodesic map created by
intersecting the normal fan of $Q$ with $S^3$.

Our construction will have the following properties, which are already
satisfied by the $Q^+$ and $Q^-$ defined in
Section~\ref{sub:28vertices}:

\begin{properties}\label{properties}
$Q^+$ and $Q^-$ contain the origin in their interior, and
the set $S^+$ of vertices of $Q^+$ (respectively, $S^-$ of $Q^-$) which project
to vertices of $\pi_{34}(Q^+)$ (resp. of $\pi_{12}(Q^-)$) consists of eight vertices, all lying
in the plane $P_{34}$ (resp. $P_{12}$).
\end{properties}

Note that these properties imply each polytope has at least eleven
vertices; e.g. $Q^+$ has eight vertices in $P_{34}$, and at least
three vertices out of $P_{34}$ so that their convex hull after
projection on $P_{12}$ contains the origin.

\begin{lemma}\label{lemma:cycle}
  Let $Q$ be a polytope in $\real^d$ with normal geodesic map $\gauss(Q)$ in $S^{d-1}$.
  Let $\Pi$ an affine subspace containing the origin and let $\pi(Q)$ be the orthogonal projection of $Q$ on $\Pi$.

  Then, every full-dimensional cell of $\gauss(Q)$ whose closure intersects $\Pi\cap S^{d-1}$
  is the normal cell of a vertex $v$ such that $\pi(v)$ is a vertex of  $\pi(Q)$.
\end{lemma}

\begin{proof}
  By definition, a vector $\phi$ is in a (closed) cell of $\gauss(Q)$
  corresponding to some vertex $v$ of $Q$ if and only if there is a
  supporting hyperplane of $Q$  perpendicular to $\phi$ and containing
  $v$. If $\phi$ is in $\Pi$, then the intersection of the
  same hyperplane with $\pi(Q)$ is $\{\pi(v)\}$; and so $\pi(v)$ is a
  vertex of $\pi(Q)$.
\end{proof}

\begin{corollary}
Let $Q^+$ be in the conditions of Properties~\ref{properties}.

  Let $D^+$ be
  the union of the closed facets in $\gauss(Q^+)$ corresponding to $S^+$. Then
  $D^+$ contains $C_{34}$ in its interior, and is partitioned by
  $\gauss(Q^+)$ into eight facets arranged in consecutive circular
  ordering around $C_{34}$, each cell only touching the preceding
  and successive one in the ordering. By choosing every second cell in
  the ordering, we get a family $F^+$ of four facets that share no
  vertices. There are no vertices of $\gauss(Q^+)$ in the
  interior of $D^+$.
\end{corollary}

\begin{figure}[ht]
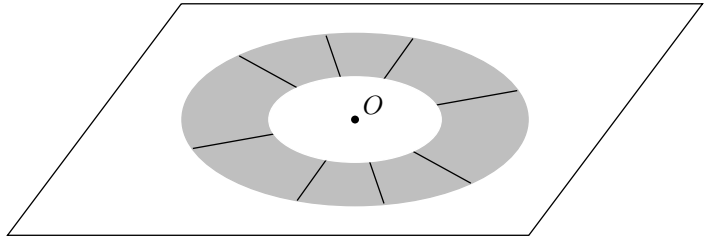

\begin{center}
  \psset{unit=0.77cm,arrows=-,shortput=nab,linewidth=0.5pt,arrowsize=2pt 5,labelsep=3.5pt}
  \pspicture(-6,-2)(6,2)
  \psline(-3,2)(6,2)(3,-2)(-6,-2)(-3,2)
  \psellipse[linecolor=lightgray,fillstyle=solid,fillcolor=lightgray](0,0)(3,1.5)
  \psline(1,1.4)(-1,-1.4)
  \psline(2.8,0.5)(-2.8,-0.5)
  \psline(2,-1.1)(-2,1.1)
  \psline(-0.5,1.45)(0.5,-1.45)
  \psellipse[linecolor=white,fillstyle=solid,fillcolor=white](0,0)(1.5,0.75)
  \psdot(0,0)
  \uput[40](0,0){$O$}
  \endpspicture
  \caption{Projection on $P_{34}$ of $D^+$ and the eight cells in $D^+$.}
  \label{fig:p34}
\end{center}
\end{figure}

\begin{proof}
$C_{34}$ intersects the eight cells of $\gauss(Q^+)$ corresponding to $S^+$, and only those,
by Lemma~\ref{lemma:cycle}.

More precisely, the normal fan (or the Gauss map) of a projection of the orthogonal projection
of a polytope $Q$ to a hyperplane (or lower dimensional linear subspace) $\Pi$
equals the intersection of the original Gauss map with $\Pi$ (in formula, and with
the notation of Lemma~\ref{lemma:cycle}, $\gauss(\pi(Q))=\gauss(Q)\cap \Pi$).
In the conditions of Properties~\ref{properties} we have that $\pi_{34}(Q^+)$ is an
octagon, so its Gauss map is a cycle of eight arcs along the circle $C_{34}$, as in Figure~\ref{fig:p34}.
This implies, as stated, that   $C_{34}$ is covered by the cells of $\gauss(Q^+)$ corresponding
to the eight vertices in $S^+$.

Since the wall between two consecutive cells in the cycle
is perpendicular to~$C_{34}$ (because the edge joining the corresponding two vertices of $Q^+$
lies in~$P_{34}$) the only way in which two non-consecutive cells in the cycle could touch each other
would be in the orthogonal complement of $C_{34}$, that is, in $C_{12}$. To finish the proof
we thus only need to show that $C_{12}$ does not touch any of the eight cells in the cycle.
This is so because $\pi_{12}(S^+)$ is just the origin,  which is in the interior of $\pi_{12}(Q^+)$, and again by Lemma~\ref{lemma:cycle}. For the same reason no vertices of $\gauss(Q^+)$ can be in the interior of $D^+$, since these would mean that (at least) four of the eight cells meet at that vertex, and no more than two of them meet.
\end{proof}

Symmetrically, if $D^-$ is the closure of the union of the facets corresponding to $S^-$ in
$\gauss(Q^-)$, then $D^-$ contains $C_{12}$, $D^-$ is partitioned into
eight facets arranged in consecutive circular ordering around $C_{12}$,
each cell only touching the preceding and successive one in the
ordering. By choosing every second cell in the ordering, we get a
family of $F^-$ of four facets with no common vertex; and no vertex of
$\gauss(Q^-)$ is in the interior of $D^-$.

Recall that the reduced incidence pattern contains only facets of a
map that contain a vertex of the other map. If two such facets have a
vertex in common, they must both have in the reduced incidence pattern
an arrow to a same facet of the other map, which makes it more
difficult to get the pattern we want. Thus, the fact that the four
facets in $F^+$ (resp. $F^-$) have no vertices in common makes them
good candidates for the reduced incidence pattern.

\begin{lemma}~\label{lemma:def}
  Any polytopes $Q^+$, $Q^-$ with Properties~\ref{properties} can
  be transformed, by independent variable scalings in $Q^+$ and $Q^-$,
  so that all vertices of $\gauss(Q^+)$ are in $D^-$ and all
  vertices of $\gauss(Q^-)$ are in $D^+$.
\end{lemma}

\begin{proof}
  We multiply by an arbitrary large number the third and fourth
  coordinates of vertices of $Q^+$, and the first and second
  coordinates of $Q^-$. This does not change the combinatorial
  properties of $\gauss(Q^+)$ or $\gauss(Q^-)$, or their positions with respect
  to $C_{12}$ and $C_{34}$,
  but deforms them by
  bringing every vertex of~$\gauss(Q^+)$ arbitrarily close to
  $C_{12}$, which is in the interior of $D^-$, and every vertex of~$\gauss(Q^-)$
  arbitrarily close to~$C_{34}$, which is in the
  interior of~$D^+$.
\end{proof}

By Lemma~\ref{lemma:def}, we can assume that all vertices of
$\gauss(Q^+)$ and $\gauss(Q^-)$ are in $D^-$ and $D^+$
respectively. The only remaining condition we need is that the
vertices of facets in $F^+$ should be in the proper facets of $F^-$
and vice-versa, so as to have the reduced incidence pattern we
seek. Recall that $D^+$ is composed of facets corresponding to
vertices in $S\subset P_{34}$, and so the partition of $D^+$ only depends on
the third and fourth coordinates. Therefore, the cell of $D^+$ that a
vertex of~$F^-$ is in is determined only by the third and fourth
coordinates of the vertex, more precisely their sign and their ratio
(see Figure~\ref{fig:p34}). Similarly, the cell of $D^-$ that a vertex
of $F^+$ is in is determined only by the first and second coordinates
of the vertex.

These properties motivate the use of diagrams presented in Section~5
of~\cite{Santos:Hirsch-counter} for representing the geodesic maps
$\gauss(Q^+)$ and $\gauss(Q^-)$ on a flat torus. The horizontal
coordinate of some vertex of a geodesic map is the circular angle defined by
its projection on the plane $\real^2\times\{(0,0)\}$; and the vertical
coordinate is the angle defined by its projection on the plane
$\{(0,0)\}\times\real^2$. These coordinates can also be thought of as
Hopf coordinates, as presented in more detail
in Section~\ref{sub:HopfCoordinates}.

\begin{figure}[ht]
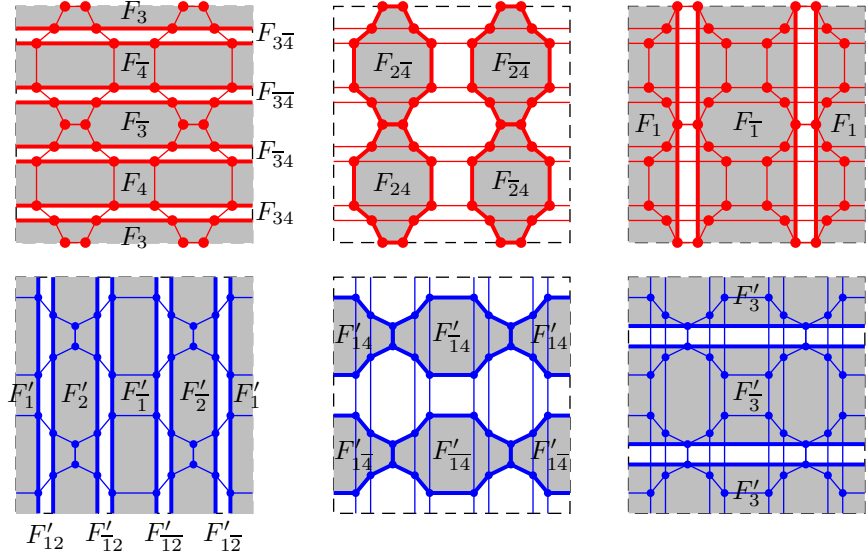

\begin{center}
  \begin{tabular}{ccc}
    \begin{minipage}{3.8cm}
      \include{Figures/newpxinc1}
     \end{minipage}&
     \begin{minipage}{3.5cm}
        \include{Figures/newpxinc1b}
     \end{minipage}&
     \begin{minipage}{3.5cm}
         \include{Figures/newpxinc1c}
     \end{minipage}\\
    \begin{minipage}{3.8cm}
      \include{Figures/newpxinc2}
     \end{minipage}&
     \begin{minipage}{3.5cm}
        \include{Figures/newpxinc2b}
     \end{minipage}&
     \begin{minipage}{3.5cm}
         \include{Figures/newpxinc2c}
     \end{minipage}\\
  \end{tabular}
  \caption{Diagram for the geodesic normal maps of polytopes $Q^+$ and $Q^-$
    of Corollary~\ref{coro:28vertices}, with 14 facets each, projected to the standard torus}
  \label{fig:Q28}
\end{center}
\end{figure}

We represent in Figure~\ref{fig:Q28} the diagrams corresponding to the
$Q^+$ (top row) and $Q^-$ (bottom row) of Section~\ref{sub:28vertices}, with 14 vertices each.
The first image in each row shows the eight facets in $D^+$ and $D^-$
respectively, with the four shaded ones being the ones in the reduced incidence diagram.

The shaded regions in the other two drawings in each row show the other six facets in $\gauss(Q^+)$
and $\gauss(Q^-)$, respectively. Thick lines in all the drawings show the separation between the facets displayed in that part.

\begin{figure}[ht]
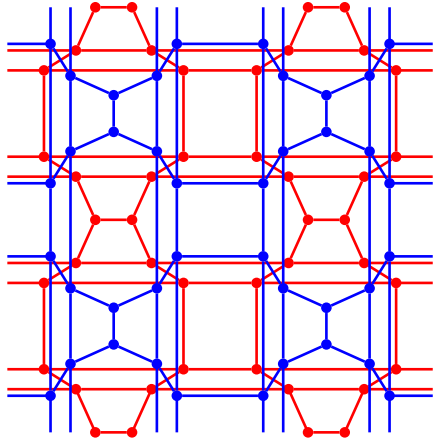

\begin{center}
  \include{Figures/newpxinc3}
  \caption{Superposition of the two geodesic maps of Corollary~\ref{coro:28vertices}}
  \label{fig:Q28-2}
\end{center}
\end{figure}

By Lemma~\ref{lemma:def}, we can assume all vertices of $\gauss(Q^+)$ and
$\gauss(Q^-)$ are in $D^-$ and $D^+$ respectively.
That is to say, shaded regions in the first image of each row are the areas where vertices of the other map are allowed. Figure~\ref{fig:Q28-2} shows that this is indeed where all the vertices lie.

In order to find smaller examples, we reduce the
number of facets in both geodesic maps. Since, by
Properties~\ref{properties}, the number of facets in $D^+$ and $D^-$
is fixed at eight, we do this by merging facets not in $D^+$ or
$D^-$. More precisely, we merge:
\begin{itemize}
\item Facets $F_{\overline 1}$ and $F_{\overline 24}$ of $Q^+$, as well as facets $F_{1}$ and $F_{\overline 2 \overline 4}$.
\item Facets $F'_{14}$ and $F'_{\overline 3}$ of $Q^-$.
\end{itemize}
 In order to keep the vertices of each geodesic map inside the appropriate facets of the other,
 the way we did this was moving little by little the vertices of
$Q^+$ and $Q^-$, using the diagram to check that the properties are
always satisfied, until we eventually got the diagram in
Figure~\ref{fig:p21}, corresponding to the following prismatoid:

\begin{figure}[ht]
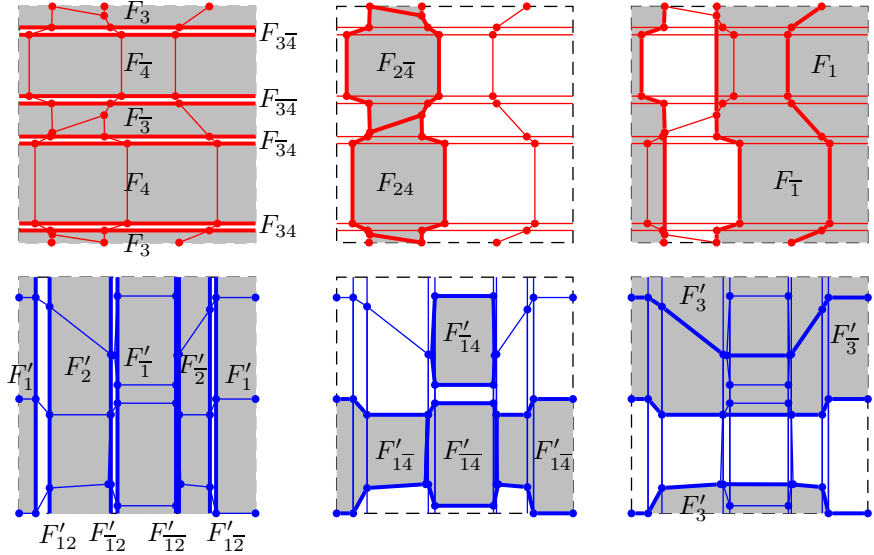

\begin{center}
  \begin{tabular}{ccc}
    \begin{minipage}{3.8cm}
      \include{Figures/s20vxinc1}
     \end{minipage}&
     \begin{minipage}{3.5cm}
        \include{Figures/s20vxinc1b}
     \end{minipage}&
     \begin{minipage}{3.5cm}
         \include{Figures/s20vxinc1c}
     \end{minipage}\\
    \begin{minipage}{3.8cm}
      \include{Figures/s20vxinc2}
     \end{minipage}&
     \begin{minipage}{3.5cm}
        \include{Figures/s20vxinc2b}
     \end{minipage}&
     \begin{minipage}{3.5cm}
         \include{Figures/s20vxinc2c}
     \end{minipage}\\
  \end{tabular}
  \caption{Diagram for the geodesic normal maps of the smaller polytopes $Q^+$ and $Q^-$,
    with 12 and 13 facets respectively, projected to the standard torus}
  \label{fig:p21}
\end{center}
\end{figure}

\begin{figure}[ht]
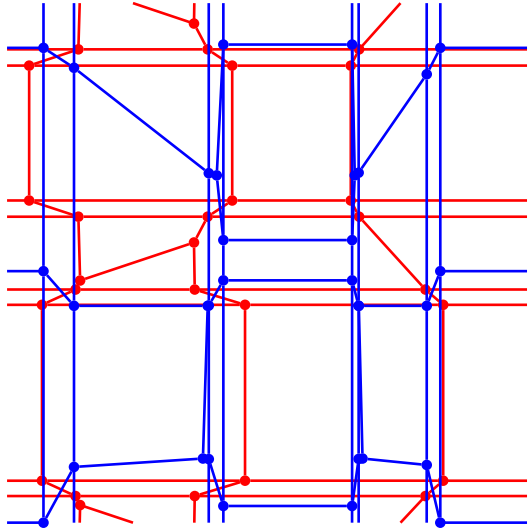

\begin{center}
  \include{Figures/s20vxinc3}
  \caption{Superposition of the two geodesic maps of the smaller polytopes $Q^+$ and $Q^-$}
  \label{fig:p21-2}
\end{center}
\end{figure}

{\footnotesize
\[
Q_{25}:= \operatorname{conv}
\left\{
\bordermatrix{
&x_1&x_2&x_3&x_4&x_5\cr
&   0&    0& 20&    -4 & 1 \cr
&   0&    0&- 20&    -4 & 1 \cr
&   0&    0& 21&    -7 & 1 \cr
&   0&    0&- 21&    -7 & 1 \cr
&   0&    0& 16&   -15 & 1 \cr
&   0&    0&- 16&   -15 & 1 \cr
&   0&    0&     0& 32 & 1 \cr
&   0&    0&     0&- 32 & 1 \cr
&\frac{ 3}{50}&\frac{-1}{25}&    0& -30 & 1 \cr
&\frac{-3}{50}&\frac{-1}{25}&    0&  30 & 1 \cr
&\frac{ 3}{1000}&\frac{7}{1000}&    0&\frac{-318}{10} & 1 \cr
&\frac{-3}{1000}&\frac{7}{1000}&    0&\frac{ 318}{10} & 1 \cr
}\quad
\bordermatrix{
&x_1&x_2&x_3&x_4&x_5\cr
&    60&     0&     0&     0 & -1 \cr
&     8&   -30&     0&     0 & -1 \cr
&     0&   -33&     0&     0 & -1 \cr
&    -2&   -32&     0&     0 & -1 \cr
&   -55&     0&     0&     0 & -1 \cr
&   -34&    36&     0&     0 & -1 \cr
&     0&    76&     0&     0 & -1 \cr
&    44&    34&     0&     0 & -1 \cr
&   -20&     0& \frac{1}{5}& \frac{-1}{5} & -1 \cr
&\frac{2999}{50}&  0& \frac{-3}{25}& \frac{-1}{5} & -1 \cr
&\frac{299999}{5000}&  0&  0& \frac{1}{100} & -1 \cr
&\frac{-549}{10}&  0& \frac{1}{5000}& \frac{1}{800} & -1 \cr
&  -54&  0& \frac{1}{500}& \frac{-1}{80} & -1 \cr
}
\right\}.
\]
} Note the differences of scale between the coefficients for $x_1$ and
$x_2$ on one side, and $x_3$ and $x_4$ on the other side. It ensures
as per Lemma~\ref{lemma:def} that vertices of~$\gauss(Q^+)$
and~$\gauss(Q^-)$ are inside $D^-$ and $D^+$ respectively.

\begin{theorem}
\label{thm:q25}
The $5$-dimensional prismatoid with $25$ vertices $Q_{25}$ described above has width six. Hence, there is a $20$-dimensional polytope with $40$ facets and diameter at least $21$.
\end{theorem}

We have seen that Properties~\ref{properties} demand that polytopes
$Q^+$ and $Q^-$ have at least eleven vertices each; they have twelve
and thirteen vertices in this construction, which is therefore close
to optimal. In fact, we believe there is no smaller $5$-prismatoid of
width $6$ with these properties.

\subsection{An explicit non-Hirsch polytope} \label{sub:Explicit}

The possibility of building and checking an explicit non-Hirsch
polytope was left as an open question
in~\cite{Santos:Hirsch-counter}. The non-Hirsch polytope of
dimension $43$ whose
existence was shown in that paper was estimated to have a
trillion vertices, which would be hard to build, and still harder to
check.  The $5$-prismatoids of width $6$ we presented in
Sections~\ref{sub:4cube} and~\ref{sub:25vertices} lead
to non-Hirsch polytopes of a much more reasonable size. It was
nevertheless necessary to find appropriate methods to build and check
them. To see how we do it, let us briefly sketch the proof of
Theorem~\ref{thm:dstep-prismatoid} contained in~\cite{Santos:Hirsch-counter}.

Let $Q_d$
be a $d$-prismatoid of dimension $d$, with $n$ vertices and width $l$, such that
$n>2d$ and $l>d$. Since $n>2d$, at least one of the two base facets, let us call it $F$,
 of $Q_d$ is not a simplex. Choose any vertex $v$ in \emph{the other} base
facet. Execute a \emph{one-point-suspension} of the prismatoid over
$v$, which consists in embedding the prismatoid in a hyperplane in
$\real^{d+1}$, and replacing $v$ by two vertices $v^+$ and $v^-$ away
from that hyperplane, such that $v$ lies in the segment from $v^+$
to $v^-$. Then perturb slightly the vertices of $F$ (the facet that was not a simplex)
away from the hyperplane.
We have the following (see details in~\cite{Santos:Hirsch-counter}[Theorem 2.6]):

\begin{lemma}
\label{lemma:suspend-and-perturb}
If the perturbation of the vertices of $F$ is done appropriately then the new polytope $Q_{d+1}$
is a $(d+1)$-prismatoid with $n+1$ vertices and width (at least) $l+1$.
\end{lemma}

Iterating this
operation $n-2d$ times successively yields a non-Hirsch prismatoid
$P_{2n-2d}$, as stated in Theorem~\ref{thm:dstep-prismatoid}.

So, to get a (dual) non-Hirsch polytope from our $5$-prismatoid with $25$
vertices we need to apply the suspension-plus-perturbation operation $15$
times. On every iteration, the perturbation needs to be done carefully and be
sufficiently small in order for the width to be augmented by one. Instead of
\emph{a priori} computing the right perturbations, we chose to check
the width of every intermediate prismatoid, changing to a smaller perturbation
when needed. For this, we
computed the facets of the prismatoid, then the adjacency list of the
facets, and computed the shortest path between the base facets. The
facets were computed efficiently using \texttt{lrs}~\cite{Avis:lrs}, which is
particularly well-suited to enumerating the facets of
nearly-simplicial polytopes. The hardest task proved to be computing
the adjacency list, considering some prismatoids we computed had tens
of thousands of facets.

The program \texttt{lrs} can be asked to output a triangulation of the
boundary of a polytope given by its vertices, and also provides the
list of vertices of each of the simplices of the triangulation. As the
prismatoids we computed are nearly simplicial (and closer to simplicial
as we went up the tower of suspensions), we were able to find
easily simplices belonging to a same facet, and thus obtain the list
of vertices incident to each facet.

From the list of vertices incident to each facet, testing whether two
facets were adjacent was done by counting the number of vertices
they had in common. In a $d$-dimensional prismatoid, adjacent facets
would have at least $d-1$ vertices in common. If the lists of vertices
incident to facets are encoded into integers, each set bit
representing an incident vertex, the list of common vertices is
computed by a simple AND operation. The number of set bits in the
result is the number of common vertices. As our prismatoids were not
simplicial, but nearly-simplicial, it was possible for some
non-adjacent facets to have $d-1$ vertices in common. In order to find
adjacent facets only, we first computed for each facet $F$ the
collection of all facets that had $d-1$ vertices in common with~$F$,
and stored for each of them the list of common vertices. The facets
adjacent with $F$ were then obtained by finding out the facets in the
collection whose list of common vertices with $F$ were maximal.

Using these methods, it was possible to build from the $5$-prismatoid
described in Section~\ref{sub:25vertices} a $20$-prismatoid of $40$
vertices and width $21$. Its vertices are the rows of the matrix in Table~\ref{table:poly20dim21}.
%

\begin{sidewaystable}
{\tiny
$\left(
\begin{tabular}{cccccccccccccccccccc}
      1    &    0    &    0    &    20    &    -4    &    0    &    0    &    0    &    0    &    0    &    0    &    0    &    0    &    0    &    0    &    0    &    0    &    0    &    0    &    0  \\
      1    &    0    &    0    &    -20    &    -4    &    0    &    0    &    0    &    0    &    0    &    0    &    0    &    0    &    0    &    0    &    0    &    0    &    0    &    0    &    1 \\
      1    &    0    &    0    &    21    &    -7    &    0    &    0    &    0    &    0    &    0    &    0    &    0    &    0    &    0    &    0    &    1    &    1    &    0    &    0    &    0\\
      1    &    0    &    0    &    -21    &    -7    &    0    &    0    &    0    &    0    &    0    &    0    &    0    &    0    &    0    &    0    &    0    &    0    &    0    &    0    &    0 \\
      1    &    0    &    0    &    16    &    -15    &    0    &    0    &    0    &    0    &    0    &    0    &    0    &    0    &    0    &    0    &    0    &    1    &    1    &    0    &    0 \\
      1    &    0    &    0    &    -16    &    -15    &    0    &    0    &    0    &    0    &    0    &    0    &    0    &    0    &    0    &    0    &    0    &    0    &    0    &    1    &    1 \\
      1    &    0    &    0    &    0    &    32    &    0    &    0    &    0    &    0    &    0    &    0    &    0    &    0    &    0    &    0    &    0    &    0    &    0    &    0    &    0 \\
\smallskip
      1    &    0    &    0    &    0    &    -32    &    0    &    0    &    0    &    0    &    0    &    0    &    0    &    0    &    0    &    0    &    0    &    0    &    1    &    1    &    0 \\
\smallskip
      1    &    $\frac{3}{50}$   &     -$\frac{1}{25}$   &    0    &    -30    &    0    &    0    &    0    &    0    &    0    &    0    &    0    &    0    &    0    &    0    &    0    &    0    &    0    &    0    &    0\\
\smallskip
      1    &     -$\frac{3}{50}$   &     -$\frac{1}{25}$   &    0    &    30    &    0    &    0    &    0    &    0    &    0    &    0    &    0    &    0    &    0    &    0    &    0    &    0    &    0    &    0    &    0  \\
\smallskip
      1    &     -$\frac{3}{1000}$   &     $\frac{7}{1000}$   &    0    &     $\frac{318}{10}$   &    0    &    0    &    0    &    0    &    0    &    0    &    0    &    0    &    0    &    1    &    0    &    0    &    0    &    0    &    0  \\
\smallskip
      1    &     $\frac{3}{1000}$   &     $\frac{7}{1000}$   &    0    &     -$\frac{318}{10}$   &    $10^7$    &    $10^7$    &    $10^7$    &    $10^{10}$    &    $10^{11}$    &    $10^{11}$    &    $10^{11}$    &    $10^{11}$    &    1    &    0    &    0    &    0    &    0    &    0    &    0   \\
\smallskip
      1    &     $\frac{3}{1000}$   &     $\frac{7}{1000}$   &    0    &     -$\frac{318}{10}$   &    -$10^7$    &    0    &    0    &    0    &    0    &    0    &    0    &    0    &    1    &    0    &    0    &    0    &    0    &    0    &    0  \\
\smallskip
      1    &     $\frac{3}{1000}$   &     $\frac{7}{1000}$   &    0    &     -$\frac{318}{10}$   &    $10^7$    &    -$10^7$    &    0    &    0    &    0    &    0    &    0    &    0    &    1    &    0    &    0    &    0    &    0    &    0    &    0  \\
\smallskip
      1    &     $\frac{3}{1000}$   &     $\frac{7}{1000}$   &    0    &     -$\frac{318}{10}$   &    $10^7$    &    $10^7$    &    -$10^7$    &    0    &    0    &    0    &    0    &    0    &    1    &    0    &    0    &    0    &    0    &    0    &    0  \\
\smallskip
      1    &     $\frac{3}{1000}$   &     $\frac{7}{1000}$   &    0    &     -$\frac{318}{10}$   &    $10^7$    &    $10^7$    &    $10^7$    &    -$10^{10}$    &    0    &    0    &    0    &    0    &    1    &    0    &    0    &    0    &    0    &    0    &    0  \\
\smallskip
      1    &     $\frac{3}{1000}$   &     $\frac{7}{1000}$   &    0    &     -$\frac{318}{10}$   &    $10^7$    &    $10^7$    &    $10^7$    &    $10^{10}$    &    -$10^{11}$    &    0    &    0    &    0    &    1    &    0    &    0    &    0    &    0    &    0    &    0  \\
\smallskip
      1    &     $\frac{3}{1000}$   &     $\frac{7}{1000}$   &    0    &     -$\frac{318}{10}$   &    $10^7$    &    $10^7$    &    $10^7$    &    $10^{10}$    &    $10^{11}$    &    -$10^{11}$    &    0    &    0    &    1    &    0    &    0    &    0    &    0    &    0    &    0  \\
\smallskip
      1    &     $\frac{3}{1000}$   &     $\frac{7}{1000}$   &    0    &     -$\frac{318}{10}$   &    $10^7$    &    $10^7$    &    $10^7$    &    $10^{10}$    &    $10^{11}$    &    $10^{11}$    &    -$10^{11}$    &    0    &    1    &    0    &    0    &    0    &    0    &    0    &    0  \\
\smallskip
      1    &     $\frac{3}{1000}$   &     $\frac{7}{1000}$   &    0    &     -$\frac{318}{10}$   &    $10^7$    &    $10^7$    &    $10^7$    &    $10^{10}$    &    $10^{11}$    &    $10^{11}$    &    $10^{11}$    &    -$10^{11}$    &    1    &    0    &    0    &    0    &    0    &    0    &    0  \\
      -1    &    60    &    0    &    0    &    0    &    0    &    0    &    1    &    1    &    0    &    0    &    0    &    0    &    0    &    0    &    0    &    0    &    0    &    0    &    0  \\
      -1    &    8    &    -30    &    0    &    0    &    0    &    0    &    0    &    1    &    1    &    0    &    0    &    0    &    0    &    0    &    0    &    0    &    0    &    0    &    0  \\
      -1    &    0    &    -33    &    0    &    0    &    0    &    0    &    0    &    0    &    1    &    1    &    0    &    0    &    0    &    0    &    0    &    0    &    0    &    0    &    0  \\
      -1    &    -2    &    -32    &    0    &    0    &    0    &    0    &    0    &    0    &    0    &    1    &    1    &    0    &    0    &    0    &    0    &    0    &    0    &    0    &    0  \\
      -1    &    -55    &    0    &    0    &    0    &    0    &    0    &    0    &    0    &    0    &    0    &    1    &    1    &    0    &    0    &    0    &    0    &    0    &    0    &    0  \\
      -1    &    -34    &    36    &    0    &    0    &    0    &    0    &    0    &    0    &    0    &    0    &    0    &    1    &    0    &    0    &    0    &    0    &    0    &    0    &    0  \\
      -1    &    0    &    76    &    0    &    0    &    0    &    0    &    0    &    0    &    0    &    0    &    0    &    0    &    0    &    0    &    0    &    0    &    0    &    0    &    0  \\
\smallskip
      -1    &    44    &    34    &    0    &    0    &    0    &    0    &    0    &    0    &    0    &    0    &    0    &    0    &    0    &    0    &    0    &    0    &    0    &    0    &    0  \\
\smallskip
      -1    &    -20    &    0    &     $\frac{1}{5}$   &     -$\frac{1}{5}$   &    0    &    0    &    0    &    0    &    0    &    0    &    0    &    0    &    0    &    0    &    0    &    0    &    0    &    0    &    0  \\
\smallskip
      -1    &     $\frac{2999}{50}$   &    0    &     -$\frac{3}{25}$    &    -$\frac{1}{5}$   &    0    &    0    &    1    &    0    &    0    &    0    &    0    &    0    &    0    &    0    &    0    &    0    &    0    &    0    &    0   \\
\smallskip
      -1    &     $\frac{299999}{5000}$    &    0    &    0    &     $\frac{1}{100}$   &    0    &    1    &    0    &    0    &    0    &    0    &    0    &    0    &    0    &    0    &    0    &    0    &    0    &    0    &    0   \\
\smallskip
      -1    &    -$\frac{549}{10}$   &    0    &     $\frac{1}{5000}$   &     $\frac{1}{800}$   &    1    &    0    &    0    &    0    &    0    &    0    &    0    &    0    &    0    &    0    &    0    &    0    &    0    &    0    &    0   \\
\smallskip
      -1    &    -54    &    0    &     $\frac{1}{500}$   &     -$\frac{1}{80}$   &    0    &    0    &    0    &    0    &    0    &    0    &    0    &    0    &    $10^5$    &    $10^7$    &    $10^7$    &    $10^7$    &    $10^8$    &    $10^8$    &    $10^9$   \\
\smallskip
      -1    &    -54    &    0    &     $\frac{1}{500}$   &     -$\frac{1}{80}$   &    0    &    0    &    0    &    0    &    0    &    0    &    0    &    0    &    -$10^5$    &    0    &    0    &    0    &    0    &    0    &    0   \\
\smallskip
      -1    &    -54    &    0    &     $\frac{1}{500}$   &     -$\frac{1}{80}$   &    0    &    0    &    0    &    0    &    0    &    0    &    0    &    0    &    $10^5$    &    -$10^7$    &    0    &    0    &    0    &    0    &    0   \\
\smallskip
      -1    &    -54    &    0    &     $\frac{1}{500}$   &     -$\frac{1}{80}$   &    0    &    0    &    0    &    0    &    0    &    0    &    0    &    0    &    $10^5$    &    $10^7$    &    -$10^7$    &    0    &    0    &    0    &    0   \\
\smallskip
      -1    &    -54    &    0    &     $\frac{1}{500}$   &     -$\frac{1}{80}$   &    0    &    0    &    0    &    0    &    0    &    0    &    0    &    0    &    $10^5$    &    $10^7$    &    $10^7$    &    -$10^7$    &    0    &    0    &    0   \\
\smallskip
      -1    &    -54    &    0    &     $\frac{1}{500}$   &     -$\frac{1}{80}$   &    0    &    0    &    0    &    0    &    0    &    0    &    0    &    0    &    $10^5$    &    $10^7$    &    $10^7$    &    $10^7$    &    -$10^8$    &    0    &    0   \\
\smallskip
      -1    &    -54    &    0    &     $\frac{1}{500}$   &     -$\frac{1}{80}$   &    0    &    0    &    0    &    0    &    0    &    0    &    0    &    0    &    $10^5$    &    $10^7$    &    $10^7$    &    $10^7$    &    $10^8$    &    -$10^8$    &    0   \\
\smallskip
      -1    &    -54    &    0    &     $\frac{1}{500}$   &     -$\frac{1}{80}$   &    0    &    0    &    0    &    0    &    0    &    0    &    0    &    0    &    $10^5$    &    $10^7$    &    $10^7$    &    $10^7$    &    $10^8$    &    $10^8$    &    -$10^9$  
\end{tabular}
\right)$
}
\caption{The vertex coordinates of a $20$-dimensional (dual) non-Hirsch polytope}
\label{table:poly20dim21}
\end{sidewaystable}

Its $36,425$ facets were computed by
\texttt{lrs} in $33$ seconds, and the $600$ million pairs of
facets were tested in $105$ seconds.

It was also possible to build from
the $5$-prismatoid described in Section~\ref{sub:28vertices} a
$23$-prismatoid of $46$ vertices and width $24$. Its $73,224$ facets
were computed by \texttt{lrs} in $49$ seconds, and the $2.6$ billion
pairs of facets were tested in $445$ seconds. Computations were done
on a laptop computer with a 2.5 GHz processor. Files describing the
  vertices, facets and adjacency lists of these polytopes are
  available on the web at the url
  \url{https://sites.google.com/site/christopheweibel/research/hirsch-conjecture}.

\section{Pairs of maps with large width} \label{sec:LargeWidth}

\subsection{Hopf coordinates} \label{sub:HopfCoordinates}

Before we start with the constructions let us fix some pleasant coordinates.
Let $C_{12}:=S^1\times\{(0,0)\}\subset S^3$ and $C_{34}:=\{(0,0)\}\times S^1\subset S^3$ as before.
Then every point of $C_{12}$ can be joined with every point of $C_{34}$ with a unique geodesic in~$S^3$.
Therefore, $S^3\iso C_{12}*C_{34} \iso (S^1)^{*2}$.

This motivates the \emph{Hopf coordinates}, which are in particular useful for visualization: Let $\T:=\R^2/m\Z^2$ be the \emph{flat torus of size $m$}. For an angle $\alpha\in[0,\pi/2]$, let
\[
T_\alpha:=\{(x,y,z,w)\in\RR^4\st x^2+y^2=\cos^2\alpha, z^2+w^2=\sin^2\alpha\}
\]
be a torus in $S^3$, which degenerates to $C_{12}$ and $C_{34}$ for $\alpha=0$ and $\alpha=\pi/2$, respectively.
We call $T_{\pi/4}$ the \emph{standard torus},
\[
T_{\pi/4}=\left\{(x,y,z,w)\in\R^4\ |\  x^2+y^2=\frac{1}{2}, z^2+w^2=\frac{1}{2}\right\}.
\]
We can parametrize $T_\alpha$ via the map
\begin{align*}
&f_\alpha:& \T\quad\to &\quad T_\alpha\cr
&& (a,b) \mapsto& \left(\cos\alpha\cos \frac{2\pi}{m} a,\cos\alpha\sin \frac{2\pi}{m} a,\sin\alpha\cos \frac{2\pi}{m} b,\sin\alpha\sin \frac{2\pi}{m} b\right).
\end{align*}
Note that $f_0$ is independent of the second argument, whereas $f_{\pi/2}$ is independent of the first argument.
This parametrizes $S^3$ as follows (in a non-injective way),
\[
F:[0,\pi/2]\times \T\to S^3, F(\alpha,a,b) = f_\alpha(a,b).
\]
This is in particular nice since the images $F(g_{ab})$ of the horizontal segments $g_{ab}:=[0,\pi/2]\times\{(a,b)\}$ are geodesics that connect the two points $f_0(a,b)\in C_{12}$ and $f_{\pi/2}(a,b)\in C_{34}$. In fact this parametrizes the geodesics $F(g_{ab})$ by arc length.

\subsection{A twisted product of two polygons} \label{sub:TwistedProduct}

We here introduce a certain $4$-polytope $P_{d,q}$ with vertices in $S^3$. We will later radially project its boundary faces to $S^3$ so that it becomes a geodesic map. In Section~\ref{sub:TopToBottom} we show that two copies of it, suitably rotated to one another, provide a pair of geodesic maps of large width.

The vertices of $P_{d,q}$ will all lie in the standard torus $T_{\pi/4}$ in $S^3$, which we defined in Section \ref{sub:HopfCoordinates}.
To simplify notation we write $f$ for $f_{\pi/4}$.

We fix two integers $d\geq 3$ and $q\geq 1$ (later we will need $q$ to be even) and let $m:=dq$. We define the vertex set $V_{d,q}$ of $P_{d,q}$ to be the image under $f$ of
\[
W_{d,q} :=
 \left\{\left({i},{j}\right)\in \T\ |\  i,j\in \Z_m, \ i-j =0 \pmod q\right\}
\]

\begin{figure}
\centerline{
\input{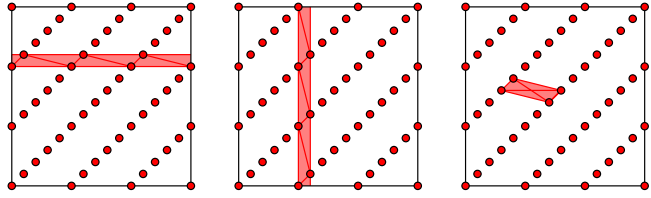}\
}
\caption{The three types of facets of $P_{3,5}$: horizontal (left), vertical (center) and diagonal (right)}
\label{fig:P35-facets}
\end{figure}

\begin{theorem}
\label{thm:facets}
$P_{d,q}$ has the following three types of facets, and only these (see Figure~\ref{fig:P35-facets}):
\begin{enumerate}
\item For each $i\in \Z_m$, the \emph{vertical facet} with vertex set the image under $f$ of
\[
 \left\{\left(i,j\right), \left({i+1},{j+1}\right) \in \T\ |\  j\in \Z_m, \ i-j =0 \pmod q\right\}.
\]
\item For each $j\in \Z_m$, the \emph{horizontal facet} with vertex set the image under $f$ of
\[
 \left\{\left(i,j\right), \left({i+1},{j+1}\right) \in \T\ |\  i\in \Z_m, \ i-j =0 \pmod q\right\}.
\]
\item For each $i,j\in \Z_m$ with $i-j =0 \pmod q$ and for each $k=1,\dots q-1$ the
\emph{diagonal tetrahedron} with vertices the image under $f$ of the following four points:
\[
\left\{\left(i,j\right),\ \left({i+1},{j+1}\right),\
\left({i+q-k},{j-k}\right),\ \left({i+q-k+1},{j-k+1}\right)\right\}.
\]
\end{enumerate}
In particular, $P_{d,q}$ has $m(n-d+2)$ facets.
\end{theorem}

Observe that if $q=1$ then $P_{d,q}$ is simply the product of two $m$-gons, in which case the horizontal and vertical facets described above are prisms and there is no diagonal tetrahedron. If $q\ge 2$ then the horizontal and vertical facets become combinatorially \emph{anti-prisms}, and diagonal tetrahedra appear connecting the two cycles of anti-prisms to one another.

\begin{proof}
We first prove that all the sets listed are indeed vertex sets of facets of~$P_{d,q}$. That there are no other facets follows from the fact that the ones listed have the property that every facet of one facet belongs also to another facet.

To show that these are facets, observe that since $V_{d,q}$ lies on the unit sphere, faces of  $P_{d,q}$ are
characterized by the same \emph{empty sphere property} that is commonly used in Delaunay triangulations: a set $X\subseteq P_{d,q}$ is the
vertex set of a face of $V_{d,q}$ if and only if there is a $2$-sphere $S_X\subset S^3$ with all points of $X$ lying on it and all points of $V_{d,q}
\setminus X$ lying on one side of it. To apply this criterion it is a bit more convenient to measure distance along the sphere $S^3$ rather than along $\RR^4$. That is, the distance between two points $x,y\in S^3$ is the angle $\arccos \langle x, y \rangle \in [0,\pi]$.

With this it is easy to show that:
\begin{enumerate}
\item The vertical facet indexed by $i$ has all its vertices in the sphere with center
\[
c= \left(\cos \frac{2\pi}{m}\left(i+\frac{1}{2}\right),\sin  \frac{2\pi}{m} \left(i+\frac{1}{2} \right), 0, 0\right) \in C_{12}
\]
and angle $\arccos \left(\frac{1}{\sqrt{2}} \cos \frac{\pi}{m} \right)$. All other vertices of $P_{d,q}$ are at a larger angle from $c$.
This is so since
\[
\left \langle c, f(a,b)\right \rangle =
\frac{1}{\sqrt{2}}\ { \cos \left [\frac{2\pi}{m} \left(a  -\left(i+\frac{1}{2} \right)\right) \right]}{}.
\]

\item The horizontal facet indexed by $j$ has all its vertices in the sphere with center
\[
c= \left(0,0,\cos \frac{2\pi}{m}\left(j+\frac{1}{2}\right),\sin  \frac{2\pi}{m} \left(j+\frac{1}{2} \right)\right) \in C_{34}
\]
and angle $\arccos \left(\frac{1}{\sqrt{2}} \cos \frac{\pi}{m} \right)$. All other vertices of $P_{d,q}$ are at a larger angle from $c$. \end{enumerate}

So, only the diagonal tetrahedra need some work. To simplify the computations we make a translation in $\T$ by the vector
$\left(- \frac{2i+1+q-k}{2}, -\frac{2j+1-k}{2}\right)$. Under this translation (which produces a rotation in $S^3$, hence it does not change the face structure of $P_{d,q}$) the vertices of $P_{d,q}$
are
\[
\left \{f(a,b)  \in T_{\pi/4} :  a+ \frac{1+q-k}{2}, b + \frac{1-k}{2}\in \Z_m, \ a-b =\frac{q}{2} \pmod q\right\}.
\]
The four vertices of the tetrahedra we want to consider become
\[
v_1:= f\left(- \frac{q-k+1}{2}, \frac{k-1}{2}\right),  \qquad v_2:= f\left(-\frac{q-k-1}{2}, \frac{k+1}{2}\right),
\]
\[
v_3:=f\left(\frac{q-k-1}{2}, - \frac{k+1}{2}\right)\quad and\quad v_4:= f\left(\frac{q-k+1}{2},  -\frac{k-1}{2}\right).
\]

For each $\alpha \in [0,\pi/2]$, let
\[
c_\alpha:=f_\alpha(0,0)=(\cos \alpha, 0, \sin \alpha,0)\in S^3.
\]

We claim that there is an $\alpha\in(0,\pi/2)$ for which $c_\alpha$ is at the same distance from the four points $v_1$, $v_2$, $v_3$,  and $v_4$, and closer to them than to the rest of vertices
of $P_{d,q}$.
This claim finishes the proof. For the claim, we compute the distance from $c_\alpha$ to any other point in the torus $T_{\pi/4}$, which is
\[
\left \langle c_\alpha, f(a,b)\right \rangle =
\frac{
\cos \alpha \cos \left (\frac{2\pi}{m} a \right)+
\sin \alpha \cos \left (\frac{2\pi}{m} b \right)
}{\sqrt{2}}\ .
\]

For every $\alpha$, $c_\alpha$ is at the same distance  from $v_1$ and $v_4$, and at the same distance from $v_2$ and $v_3$. Moreover, for $\alpha=0$ we have
\[
\left \langle c_0, v_1\right \rangle =
\frac{ \cos \left [\frac{2\pi}{m} \left(\frac{q-k+1}{2}\right) \right]}{\sqrt{2}} <
\frac{ \cos \left [\frac{2\pi}{m} \left(\frac{q-k-1}{2}\right) \right] }{\sqrt{2}}=
\left \langle c_0, v_2\right \rangle,
\]
while for $\alpha=\pi/2$
\[
\left \langle c_{\pi/2}, v_1\right \rangle =
\frac{\cos \left [\frac{2\pi}{m} \left(\frac{k-1}{2}\right) \right] }{\sqrt{2}} >
\frac{\cos \left [\frac{2\pi}{m} \left(\frac{k+1}{2}\right) \right] }{\sqrt{2}} =
\left \langle c_{\pi/2}, v_2\right \rangle.
\]
Hence, there is a value of $\alpha \in (0,\pi/2)$ for which $c_\alpha$ is at the same distance from the four points. We fix $\alpha$ to that value for the rest of the proof. We now need to show that the distance from $c_{\alpha}$ to any other vertex of $P_{d,q}$ is larger. For this, let $f(a,b)$ be a vertex. If either $\left|a\right|> \frac{q}{2}$ or $\left|b\right|>\frac{q}{2}$, let $(a',b')\in \T$ be a (typically unique) point with $a-a',b-b' =0 \pmod q$ and  $\left|a'\right|\le \frac{q}{2}$ and $\left|b'\right|\le\frac{q}{2}$.
Then~$f(a',b')$ is also a vertex of $P_{d,q}$ and, by the above expression for $\left \langle c_\alpha, f(a,b)\right \rangle$, its distance to $c_{\alpha}$ is strictly smaller than that of $f(a,b)$.

We finally deal with a vertex $f(a,b)$ of $P_{d,q}$ such that $\left|a\right|\le \frac{q}{2}$ and $\left|b\right|\le\frac{q}{2}$. This point lies along one of the helixes $\{f(a,b) : a-b= -q/2\}$ or $\{f(a,b) : a-b = q/2\}$ and, by symmetry, we may assume without loss of generality that it lies in the first one. So, our point is $f\left(t- \frac{q}{2}, t\right)$ for some $t\in \R/m \Z$. Plugging this into the expression for $\left \langle c_\alpha, f(a,b)\right \rangle$ we get that the distance from this point to $c_\alpha$ equals
\[
\phi(t) := \left \langle c_{\pi/2}, f\left(t- \frac{q}{2}, t\right)\right \rangle = \frac{
\cos \alpha \cos \left [\frac{2\pi}{m} \left(t  -\frac{q}{2}\right) \right] +
\sin \alpha \cos \left [\frac{2\pi}{m} t \right]
}{\sqrt{2}}\ .
\]

The curve $\left\{  \left( \cos \left [\frac{2\pi}{m} \left(t  -\frac{q}{2}\right)\right] ,
 \cos \left [\frac{2\pi}{m} t \right]\right) : t \in  \R/m \Z\right\} \subset \R^2$ is an ellipse and $\phi$ is a linear function on it. Hence, $\phi$ is \emph{unimodal}: it has a unique minimum and a unique maximum and it is monotone in the two arcs between them. Since $\phi$ takes the same value in the two $t$'s corresponding to the points $v_1$ and $v_2$ and since these points are consecutive among those of the helix that correspond to vertices of $P_{d,q}$, $\phi(t) - \left \langle c_0, v_1\right \rangle$ has the same sign on all those points. We saw that this sign is positive (for example) for the point $f\left(- \frac{q-k+1}{2}+q, \frac{k-1}{2}+q\right)$, so it is positive all throughout.
\end{proof}

\subsection{Two copies of the twisted product} \label{sub:TopToBottom}

In this section we put two affine copies of $P_{d,k}$ together to obtain a prismatoid with large width.
As in the previous section we fix two integers $d\geq 3$ and $q\geq 2$, $q$
being even this time, and let $m:=dq$.

Let $\alpha>0$ be small enough and put $\beta:=\pi/2 - \alpha$. We define the vertex sets $V^+:=f_\alpha(W^+)$ and $V^-:=f_\beta(W^-)$, where
\[
W^+ := \left\{\left({i},{j}\right)\in \T\ |\  i,j\in \Z_m, \ i-j =0 \pmod q\right\}
\]
and
\[
W^- := \left\{\left({i},{j}\right)\in \T\ |\  i,j\in \frac{1}{2}+\Z_m, \ i-j =\frac{q}{2} \pmod q\right\}.
\]

Let $Q^+$ and $Q^-$ be the polytopes polar to $\conv V^+$ and $\conv V^-$. Then:

\begin{theorem}
\label{thmPrismatoidWithArbitraryLargeWidth}
The polytope $Q:=\conv(Q^+\times \{1\} \cup Q^- \times \{-1\})$
is a $5$-dimensional prismatoid with $m(m-d+2)$ vertices in each base facet and of width $4+q/2$.
\end{theorem}

That $Q$ is a prismatoid is obvious and the number of vertices follows from Theorem~\ref{thm:facets}: vertices of $Q^+$ or $Q^-$ correspond to facets of their polars, which are affinely equivalent to the twisted product $P_{d,q}$ of the previous section.
To study the width, let $\gauss^+$ and $\gauss^-$ be the normal fans of $Q^+$ and $Q^-$ intersected with $S^3$, which are geodesic maps with the face lattice of $\conv V^+$ and $\conv V^-$.

The distance from the top to the bottom facet of the prismatoid is $2$ plus the distance of $V^+$ to $V^-$ along the graph, say $H$, of the common refinement of $\gauss^+$ and $\gauss^-$.
From the description of the facets of $\gauss^+$ (and hence $\gauss^-$) in Theorem~\ref{thm:facets} we can also describe all its lower dimensional faces:
\begin{itemize}
\item The $2$-faces (ridges) of $\gauss^+$ are the $m$ horizontal and $m$ vertical $d$-gons between consecutive horizontal and vertical facets, plus the $2m^2$ triangles bounding the diagonal tetrahedra.
\item The edges are the $md$ vertical and $md$ horizontal edges on those $d$-gons plus the $md$ short edges between points $f(i,j)$ and $f(i+1,j+1)$ within each helix and plus the $md(q-1)$ edges connecting consecutive helixes.
\end{itemize}

Our choice of $\alpha$ very small and $\beta=\pi/2 - \alpha$ makes all the vertices, edges, and non-horizontal ridges and facets of $\gauss^+$ be contained in the solid torus formed by the vertical facets and ridges of $\gauss^-$, and all the vertices, edges, and non-vertical ridges and facets of $\gauss^-$ be contained in the solid torus formed by the horizontal facets and ridges of $\gauss^+$.
See Figure \ref{fig:G+AndG-} and Lemma \ref{lemExplicitAlpha}.
Hence, the vertices of $H$ are:

\begin{figure}
\centerline{
\input{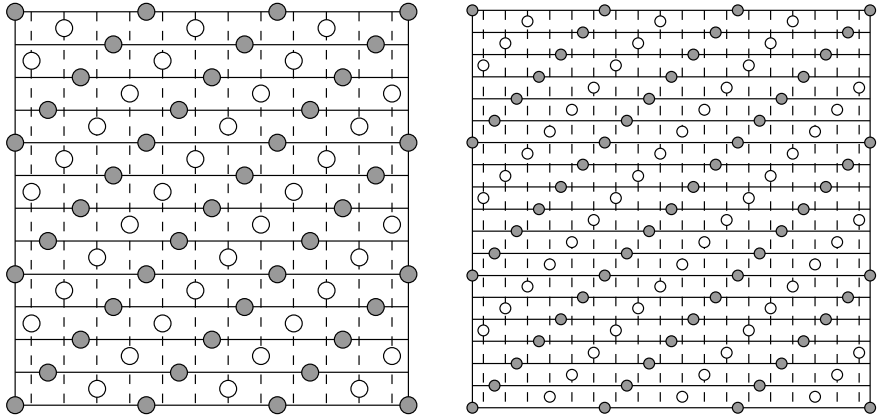}\
}
\caption{This figure shows the standard torus $T_{\pi/4}\subset S^3$ and how~$G^+$ and~$G^-$ for $d=3$ and $q=4$ (respectively $d=3$ and $q=6$) intersect it, provided $\alpha$ is small enough.
The horizontal solid lines are intersections of $T_{\pi/4}$ with the horizontal $d$-gons of $\gauss^+$ and the vertical dashed lines are intersections of $T_{\pi/4}$ with vertical $d$-gons of $\gauss^-$.
The shaded vertices are not the vertices of $\gauss^+$ but the intersections of $T_{\pi/4}$ and the geodesic segments from the midpoints of the horizontal $d$-gons of $\gauss^+$ to its vertices; analogously with the white vertices and~$\gauss^-$.
Note that $T_{\pi/4}$ separates the vertices of $\gauss^+$ from the vertices of $\gauss^-$ in~$S^3$.}
\label{fig:G+AndG-}
\end{figure}

\begin{itemize}
\item the vertices in $V^+$,
\item the vertices in $V^-$,
\item the intersections of non-vertical edges of $\gauss^+$ with vertical $d$-gons of $\gauss^-$, and
\item the intersections of non-horizontal edges of $\gauss^-$ with horizontal $d$-gons of $\gauss^+$.
\end{itemize}

\noindent
The edges of $H$ are
\begin{itemize}
\item the vertical edges of $\gauss^+$,
\item the horizontal edges of $\gauss^-$,
\item the pieces in which the non-vertical edges of $\gauss^+$ have been cut by vertical $d$-gons of $\gauss^-$,
\item the pieces in which the non-horizontal edges of $\gauss^-$ have been cut by horizontal $d$-gons of $\gauss^+$,
\item the intersections of triangles of $\gauss^+$ with vertical $d$-gons of $\gauss^-$,
\item the intersections of triangles of $\gauss^-$ with horizontal $d$-gons of $\gauss^+$,
\item the intersections of horizontal $d$-gons of $\gauss^+$ with vertical $d$-gons of $\gauss^-$.
\end{itemize}

\begin{figure}[tbh]
\centering
\input{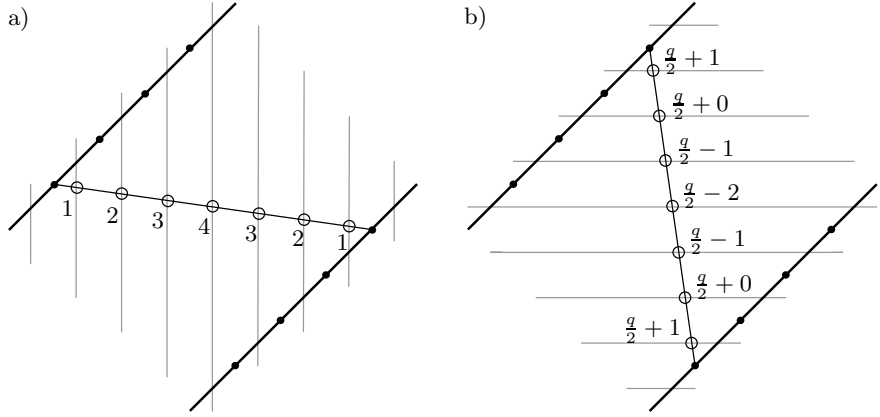}
\caption{a) Definition of $d(v)$ for intersections $v$ of the edges of $\gauss^+$ of negative slope with vertical $d$-gons of $\gauss^-$.
b) Definition of $d(v)$ for intersections $v$ of the edges of $\gauss^-$ of negative slope with horizontal $d$-gons of $\gauss^+$.}
\label{figDefinitionOfD}
\end{figure}

\noindent
We define a function $d:V(H)\to \N$ on the vertices of $H$ as follows.
\begin{itemize}
\item For $v\in V^+$, we define $d(v):=0$.
\item For the intersections $v$ of short edges of $\gauss^+$ within a helix with vertical $d$-gons of $\gauss^-$, we define $d(v)=1$.
\item For the intersections $v$ of the other edges $e$ of $\gauss^+$ with vertical $d$-gons of $\gauss^-$, we define $d(v)$ as in Figure \ref{figDefinitionOfD}a): We start at the end-points of $e$ and walk towards the midpoint of $e$. At the $i$th occurring crossing $v$ with vertical $d$-gons of $\gauss^-$ we set $d(v):=i$.
\item For $v\in V^-$, we define $d(v):=2+q/2$.
\item For the intersections $v$ of short edges of $\gauss^-$ within a helix with horizontal $d$-gons of $\gauss^+$, we define $d(v)=2+q/2$.
\item For the intersections $v$ of the other edges $e$ of $\gauss^-$ with horizontal $d$-gons of $\gauss^+$, we define $d(v)$ as in Figure \ref{figDefinitionOfD}b): We start at the end-points of $e$ and walk towards the midpoint of $e$. At the $i$th occurring crossing $v$ with horizontal $d$-gons of $\gauss^+$ we set $d(v):=2+q/2-i$.
\end{itemize}

Furthermore we add edges to $H$ between all pairs of vertices $v$ and $w$ that satisfy $d(v)=d(w)\pm 1$. We denote the resulting graph by $\wt H$.

\begin{lemma}
$d$ fulfills the properties
\begin{enumerate}
\item \label{itemProperty1} $d(V^+)=0$,
\item \label{itemProperty2} every vertex $v\in V(\wt H)\wo V^+$ has a neighbor $w$ in $\wt H$ such that $d(v)=d(w)+1$, and
\item \label{itemProperty3} no vertex $v\in V(\wt H)\wo V^+$ has a neighbor $w$ in $\wt H$ such that $d(v)\geq d(w)+2$.
\end{enumerate}
Therefore, $d$ is the distance function $\dist_{\wt H}(\underline{\ \ },V^+)$ in the graph~$\wt H$.
\end{lemma}

\begin{figure}[tbh]
\centering
\input{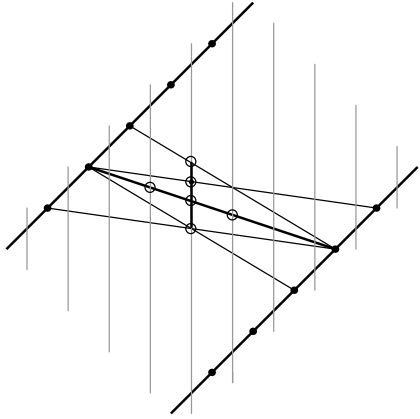}
\caption{The $6$ neighbors of $v$, two of which lie above each other.}
\label{figNeighborsOfE}
\end{figure}

\begin{proof}[Sketch of proof]
The first two properties follow immediately from the construction. For the second property let us only consider the case when $v$ is the intersection of a vertical $d$-gon $R$ of $\gauss^-$ with an edge $e$ of $\gauss^+$ with negative slope. This edge $e$ lies in $4$ ridges of $\gauss^+$, say $R_1,\ldots,R_4$. The outgoing edges from $v$ are $R\cap R_1,\ldots,R\cap R_4$, and two pieces of $e$. Therefore $e$ has $6$ neighbors, only $3$ of which have a smaller $d$-value, which is in all three cases less by $1$. Compare with Figure~\ref{figNeighborsOfE}.
For all other vertices $v$ in $\wt H$, Property \ref{itemProperty3} can be shown similarly.

The first three properties imply by induction that $d(v)=\dist_{\wt H}(v,V^+)$: Let $A_i:=\{v\st d(v)=i\}$ and $B_i:=\{v\st \dist_{\wt H}(v,V^+)=i\}$.
$A_0\supseteq B_0$ follows from Property \ref{itemProperty1} and $A_0\subseteq B_0$ follows from Property \ref{itemProperty2} and the non-negativity of $d$. Hence $A_0=B_0$.

Now suppose that $A_0=B_0,\ldots,A_i=B_i$. By Property \ref{itemProperty2}, every vertex $v\in A_{i+1}$ has a neighbor $w\in A_i=B_i$, that is, $\dist_{\wt H}(v,V^+)\leq i+1$.
Therefore $A_{i+1}\subseteq B_0\cup\ldots\cup B_{i+1}$. By disjointness of the $A_i$s, $A_{i+1}\subseteq B_{i+1}$.
To show the other inclusion, pick a new vertex $v$ in $B_{i+1}$. It has a neighbor $w$ in $B_i=A_i$, hence $d(w)=i$. By Property \ref{itemProperty3}, $d(v)\leq i+1$. Therefore $B_{i+1}\subseteq A_0\cup\ldots\cup A_{i+1}$. By disjointness of the $B_i$s, $B_{i+1}\subseteq A_{i+1}$. This finishes the induction step $A_{i+1}=B_{i+1}$.
\end{proof}

By the lemma, $d(v)=\dist_{\wt H}(v,V^+)$ for all vertices $v$ in $\wt H$.
In particular, $\dist_{\wt H}(V^+,V^-)=2+q/2$. Since $H$ is a subgraph of~$\wt H$, the distance $\dist_H(V^+,V^-)$ in $H$ is at least $2+q/2$. Since we easily  find a path from $V^+$ to $V^-$ in $H$ with that length, we have
$\dist_H(V^+,V^-)=2+q/2$.
This proves Theorem~\ref{thmPrismatoidWithArbitraryLargeWidth}. \qed

Taking $d=3$ in Theorem~\ref{thmPrismatoidWithArbitraryLargeWidth} and letting $k=q/2$ vary we get:

\begin{corollary}
\label{corPrismatoidWithArbitraryLargeWidth}
For every $k$ there is
a $5$-dimensional prismatoid with $6k(6k-1)$ vertices in each base facet and of width $4+k$.
\end{corollary}

\subsection{Explicit construction} \label{subExplicitAlpha}

The construction of $5$-prismatoids with arbitrary large width in Theorem \ref{thmPrismatoidWithArbitraryLargeWidth} is already explicit up to the choice of a sufficiently small $\alpha>0$.
The following lemma implies that $\alpha=\pi/6$ does work for all $d\geq 4$.

\begin{lemma}
\label{lemExplicitAlpha}
Suppose that $d\geq 4$.
Then any $\alpha\in (0,\pi/4)$ satisfying
\begin{equation}
\label{eqExplicitAlpha}
\tan^2 \alpha < \frac{1+\cos(2\pi/d)}{2}
\end{equation}
works for the construction of Theorem \ref{thmPrismatoidWithArbitraryLargeWidth}.
\end{lemma}

\begin{proof}
We need to show that the vertices, edges, and non-horizontal ridges and facets of $\gauss^+$ are contained in the solid torus formed by the vertical facets and ridges of $\gauss^-$, and the same with $\gauss^+$ and $\gauss^-$ interchanged.

Any vertex of $\gauss^+$ and every point on the edges and non-horizontal ridges and facets of $\gauss^+$ is a geodesic convex combination of points $P_0,\ldots,P_n\in T_\alpha$,
\[
\begin{array}{rcl}
P_i:=P_i(a,b) & := & f_\alpha(\frac{m}{2\pi}a_i,\frac{m}{2\pi}b_i) \\
              & = & (\cos\alpha\cos a_i, \cos\alpha\sin a_i, \sin\alpha\cos b_i, \sin\alpha\sin b_i),
\end{array}
\]
such that the $a_i\in [a_0,a_0+2\pi/d]$ and $b_i\in [0,2\pi]$, $0\leq i\leq n$.
A geodesic convex combination of $P_0,\ldots,P_n$ is of the form $P/||P||$, where
\[
P:=P(\lambda,a,b):=\sum_{i=0}^n\lambda_i P_i, \ \ \ \lambda_i\geq 0, \ \ \ \sum_{i=0}^n \lambda_i = 1.
\]
It is enough to prove that under the assumption \eqref{eqExplicitAlpha} all such points $P/||P||$ lie inside the solid torus
\[
T_{\leq \pi/4} := \bigcup_{0\leq x\leq \pi/4} T_x.
\]
We may assume that $a_0=0$, otherwise we use a rotation, and that $\lambda_i>0$ for all~$i$, otherwise we just omit the redundant point $P_i$ in $P_0,\ldots,P_n$.
We can explicitly compute the angle $x:=x(\lambda,a,b)$ such that $P/||P||\in T_x$ as follows:
Define the function
\[
g(\lambda,a):=\left(\sum_{i=0}^n \lambda_i\cos(a_i)\right)^2 + \left(\sum_{i=0}^n \lambda_i\sin(a_i)\right)^2.
\]
Then
\[
||P||^2 = \cos^2(\alpha) g(\lambda,a) + \sin^2(\alpha) g(\lambda,b),
\]
and
\begin{equation}
\label{eqBoundOnCosSqX}
\cos^2 x
= \frac{\cos^2(\alpha) g(\lambda,a)}{||P_\lambda||^2}
= \frac{\cos^2(\alpha) g(\lambda,a)}{\cos^2(\alpha) g(\lambda,a)+\sin^2(\alpha) g(\lambda,b)}.
\end{equation}
We have to show that \eqref{eqExplicitAlpha} implies $x\leq \pi/4$ for all $a,b,\lambda$ as above.
For this we will bound $g(\lambda,b)$ from above and $g(\lambda,a)$ from below.
The addition law for cosine yields
\[
g(\lambda,b) = \sum_{0\leq i,j\leq n} \lambda_i\lambda_j\cos(b_i-b_j).
\]
Thus,
\begin{equation}
\label{eqBoundOnGOfLambdaA}
g(\lambda,b) \leq \sum_{0\leq i,j\leq n} \lambda_i\lambda_j = \left(\sum_i \lambda_i\right)^2 = 1.
\end{equation}
Equality is attained if all $b_i$ with $\lambda_i>0$ are equal.

To get a lower bound we first note that if $d\geq 4$ then $\cos(a_i)\geq 0$ and $\sin(a_i)\geq 0$.
Thus
\[
\frac{\partial^2}{\partial a_i^2}g(\lambda,a)
=-2\lambda_i\left(\cos(a_i)\bigg(\sum_{j\neq i}\lambda_j\cos(a_j)\bigg) + \sin(a_i)\bigg(\sum_{j\neq i}\lambda_j\sin(a_j)\bigg)\right) \leq 0.
\]
Hence we may assume that $a_i\in\{0,2\pi/d\}$ for all $i$.
Thus for the lower bound on $g(\lambda,a)$ we may assume that $n=1$, $a_0=0$ and $a_1=2\pi/d$.
Using Lagrange multipliers, the minimum of $g$ is attained where
\[
\bigg(\frac{\partial}{\partial \lambda_0}-\frac{\partial}{\partial \lambda_1}\bigg) g(\lambda,a)|_{a_0=0,a_1=2\pi/d}
= (\lambda_0-\lambda_1)(1-\cos(2\pi/d))
\]
is zero, that is, at the point $\lambda_0=\lambda_1=1/2$.
Therefore,
\begin{equation}
\label{eqBoundOnGOfLambdaB}
g(\lambda,a)\geq g((1/2,1/2),(0,2\pi/d)) = (1+\cos(2\pi/d))/2.
\end{equation}
Now \eqref{eqBoundOnCosSqX}, \eqref{eqBoundOnGOfLambdaA}, and \eqref{eqBoundOnGOfLambdaB} imply
\[
\cos^2(x)\geq \frac{\cos^2(\alpha)(1+\cos(2\pi/d))/2}{\cos^2(\alpha)(1+\cos(2\pi/d))/2+\sin^2(\alpha)}.
\]
Hence we are done if the right hand side of the last inequality is larger than $\cos^2(\pi/4)=1/2$.
This condition is equivalent to the assumption~\eqref{eqExplicitAlpha}.
\end{proof}

\vskip.5cm
{\small
\noindent {Benjamin Matschke}\newline
\emph{Forschungsinstitut f\"ur Mathematik, HG 36.1}\newline
\emph{ETH Z\"urich, CH-8092 Z\"urich, Switzerland}\newline
\emph{email: }\url{benjamin.matschke@math.ethz.ch}\newline
\emph{web: }\url{http://www.math.ethz.ch/~matschkb/}
}
\vskip.5cm
{\small
\noindent {Francisco Santos}\newline
\emph{Departamento de Matem\'aticas, Estad\'istica y Computaci\'on}\newline
\emph{Universidad de Cantabria, E-39005 Santander, Spain}\newline
\emph{email: }\url{francisco.santos@unican.es}\newline
\emph{web: }\url{http://personales.unican.es/santosf/}
}
\vskip.5cm
{\small
\noindent {Christophe Weibel}\newline
\emph{Computer Science Department, 6211 Sudikoff Lab}\newline
\emph{Dartmouth College, Hanover NH 03755, United States}\newline
\emph{email: }\url{christophe.weibel@gmail.com}\newline
\emph{web: }\url{https://sites.google.com/site/christopheweibel}
}

\end{document}